\documentclass[12pt]{article}
\usepackage{multicol}
\usepackage{amsmath, amssymb, amscd, enumerate, amsfonts, latexsym, amsxtra, 
amsthm, bm}
\usepackage{indentfirst}
\newtheorem{The}{Theorem}[section]
\newtheorem{Prop}[The]{Proposition}
\newtheorem{Lem}[The]{Lemma}

\newtheorem{Def}[The]{Definition}

\begin{document}
\allowdisplaybreaks[3]
\centerline{\Large Generating sets of the Jacobson radical of} \vspace{3mm}
\centerline{\Large the hyperalgebra of $({\rm SL}_2)_r$  } \vspace{7mm}
\centerline{Yutaka Yoshii 
\footnote{ E-mail address: yutaka.yoshii.6174@vc.ibaraki.ac.jp}}  \vspace{5mm}
\centerline{College of Education,   
Ibaraki University,}
\centerline{2-1-1 Bunkyo, Mito, Ibaraki, 310-8512, Japan}
\begin{abstract}
We give   generating sets of the Jacobson radical of the 
hyperalgebra of the $r$-th Frobenius kernel of the 
algebraic group ${\rm SL}_2$ over an algebraically closed field  of characteristic 
$p>0$. This result generalizes earlier work by Wong for $r=1$ and odd $p$. 
\end{abstract}

\section{Introduction}
Let $k$ be an algebraically closed field of characteristic $p>0$.  Let $G$ be a 
simply connected and simple algebraic group over $k$ and $G_r$ the 
$r$-th Frobenius kernel of $G$. Let $\mathcal{U}_r$ be the hyperalgebra of $G_r$. 

In ring theory, knowing about the Jacobson radical is one of the important problems.  
But unfortunately, little is known about the Jacobson radical of 
$\mathcal{U}_r$. 
The exception is the simplest case $G=\textup{SL}_2$ (i.e., $G$ is of type $\textup{A}_1$) 
and $r=1$, where Wong \cite{wong83} 
gave generating sets for the Jacobson radical of $\mathcal{U}_1$ 
when $p$ is odd. 

The aim of this paper is to produce generating sets for the Jacobson 
radical of $\mathcal{U}_r$ (for any $r$ and $p$) when $G=\textup{SL}_2$. In Section 2, 
we first describe the basic properties of the hyperalgebra $\mathcal{U}_r$ for $G=\textup{SL}_2$. Then, in Section 3, we describe the properties of the 
elements $B^{({\bm \varepsilon})}({\bm a}, {\bm j})$ in the hyperalgebra $\mathcal{U}_r$  
which the author constructed in the earlier works \cite{yoshii17}, \cite{yoshii18}, and 
\cite{yoshii22}. These elements, whose construction method is motivated by Seligman's paper 
\cite{seligman03} on primitive idempotents for $\mathcal{U}_1$, 
help us to show that each of the  generating sets 
in fact generates the Jacobson radical. In particular, in \cite[\S 5]{yoshii22}, 
the author has constructed a $k$-basis of the Jacobson radical in terms of the elements. 
That fact is stated in Theorem \ref{BasisThm} in this article, which plays an important role 
in proving the main result. Then, in Section 4, we finally state the main result. 
The result is stated in Theorem \ref{MainThm}. 
In the  theorem, if $p$ is odd, we give    
generating sets consisting of $2r$ elements, which improves Wong's result (for 
$r=1$ and odd $p$). Unfortunately, the argument for odd $p$ does not work 
when $p=2$.  Nevertheless,  we give another generating set of the radical of $\mathcal{U}_r$ 
for $p=2$ instead, though the number of generators is not $2r$ but $2(2^r-1)$.

\section{Preliminaries}
Let $k$ be an algebraically closed field of characteristic $p>0$ and  $\mathbb{F}_p$ 
a finite field of $p$ elements. 
Throughout this paper, all modules for an associative 
 $k$-algebra are assumed to be finite-dimensional left 
modules. 
For a finite-dimensional (associative) $k$-algebra $R$, let ${\rm rad}R$ be the 
largest nilpotent two-sided ideal of $R$, which is called the Jacobson radical of $R$. 
For an $R$-module $M$, 
the $R$-submodule of $M$ consisting of the elements annihilated by ${\rm rad} R$ 
is denoted by ${\rm soc}_R M$ and called the  socle of $M$. 
This is  the largest semisimple $R$-submodule of $M$. For details, for example, see 
\cite[ch. I. 1]{alperinbook}. \\

From here to the end of this article, 
let $G={\rm SL}_2$ be the special linear group of degree $2$ over $k$.  Let 
$$X= \left( \begin{array}{cc}
         0 & 1 \\
         0 & 0 \end{array} \right),\ \ \ 
Y= \left( \begin{array}{cc}
         0 & 0 \\
         1 & 0 \end{array} \right),\ \ \ 
H= \left( \begin{array}{cc}
         1 & 0 \\
         0 & {-1} \end{array} \right)$$
be the standard $\mathbb{C}$-basis in the simple complex Lie algebra $\mathfrak{g}_{\mathbb{C}}=\mathfrak{sl}_2(\mathbb{C})$. Let 
$\mathcal{U}_{\mathbb{C}}$ be the universal enveloping algebra of 
 $\mathfrak{g}_{\mathbb{C}}$. Let $\mathcal{U}_{\mathbb{Z}}$ be the  subring 
 of  $\mathcal{U}_{\mathbb{C}}$  
generated by all $X^{(m)}=X^m/m!$ and  $Y^{(m)}=Y^m/m!$ with $m \in \mathbb{Z}_{\geq 0}$. 
For $z \in \{H, -H \}$, set 
$${z+c \choose m }= \dfrac{\prod_{i=0}^{m-1}(z+c-i)}{m!}$$
for $c \in \mathbb{Z}$ and $m \in \mathbb{Z}_{\geq 0}$. 
The $k$-algebra 
$\mathcal{U}_{\mathbb{Z}} \otimes_{\mathbb{Z}} k$  can be identified with  
the hyperalgebra of $G$ and is denoted by $\mathcal{U}$ or 
${\rm Dist}(G)$. We  use the same notation for the images in $\mathcal{U}$ 
of the elements 
in $\mathcal{U}_{\mathbb{Z}}$. Then $\mathcal{U}$ has 
$Y^{(m)} {H \choose n} X^{(m')}$ 
with $m, m', n \in \mathbb{Z}_{\geq 0}$ as a $k$-basis and is a graded $k$-algebra 
by regarding each basis element $Y^{(m)} {H \choose n} X^{(m')}$ as a monomial with degree 
$m'-m$. 

Consider an involutive 
 ring automorphism of $\mathcal{U}_{\mathbb{Z}}$ defined by 
$X \mapsto -Y$, $Y \mapsto -X$ and an involutive  ring antiautomorphism 
of $\mathcal{U}_{\mathbb{Z}}$ defined by 
$X \mapsto -X$, $Y \mapsto -Y$. These maps  
induce a unique  $k$-algebra 
automorphism $\mathcal{T}_1$ of $\mathcal{U}$ and a unique  $k$-algebra antiautomorphism 
$\mathcal{T}_2$ of $\mathcal{U}$ respectively satisfying 
$$\mathcal{T}_1\left( X^{(m)} \right) = (-1)^m Y^{(m)},\ \ \ 
\mathcal{T}_1\left( Y^{(m)} \right) = (-1)^m X^{(m)},
\ \ \ \mathcal{T}_1 \left( {H \choose m}\right) = {-H \choose m}$$
and 
$$\mathcal{T}_2 \left( X^{(m)} \right) = (-1)^m X^{(m)},\ \ \ 
\mathcal{T}_2 \left( Y^{(m)} \right) = (-1)^m Y^{(m)},
\ \ \ \mathcal{T}_2 \left( {H \choose m}\right) = {-H \choose m}$$
for $m \in \mathbb{Z}_{\geq 0}$. 
On the other hand, let ${\rm Fr}: \mathcal{U} \rightarrow \mathcal{U}$ be 
 the $k$-algebra endomorphism 
 which satisfies 
$${\rm Fr}\left( X^{(m)} \right) = 
\left\{ \begin{array}{cl}
         {X^{(m/p)}} & {\mbox{if $p \mid m$,}} \\
         0 & {\mbox{if $p \nmid m$}} \end{array} \right.,  
\ \ \ \ \ \ 
{\rm Fr}\left( Y^{(m)} \right) = 
\left\{ \begin{array}{cl}
         {Y^{(m/p)}} & {\mbox{if $p \mid m$,}} \\
         0 & {\mbox{if $p \nmid m$}} \end{array} \right. ,$$
and 
$${\rm Fr} \left( {H \choose m} \right) = 
\left\{ \begin{array}{cl}
         {{H \choose m/p}} & {\mbox{if $p \mid m$,}} \\
         0 & {\mbox{if $p \nmid m$}} \end{array} \right. $$
for $m \in \mathbb{Z}_{\geq 0}$.

Let $\mathcal{U}^+$ (resp. $\mathcal{U}^-$) be the $k$-subalgebra of $\mathcal{U}$ 
generated by $X^{(n)}$ (resp.  $Y^{(n)}$) with $n \in \mathbb{Z}_{\geq 0}$ 
, and let $\mathcal{U}^0$ be 
the $k$-subalgebra of $\mathcal{U}$ generated by ${H \choose n}$ with 
$n \in \mathbb{Z}_{\geq 0}$. 
For a positive integer 
$r \in \mathbb{Z}_{> 0}$, let $\mathcal{U}_r$ be the $k$-subalgebra of $\mathcal{U}$ generated by 
$X^{(n)}$ and $Y^{(n)}$ with $0 \leq n  \leq p^r-1$. 
This is a finite-dimensional $k$-algebra of dimension $p^{3r}$ which has 
$Y^{(m)} {H \choose n} X^{(m')}$ 
with $0 \leq m, m', n \leq p^r-1$ as a $k$-basis, 
and it can be identified with the hyperalgebra of the $r$-th Frobenius kernel 
$G_r$ of $G$. 
Set $\mathcal{U}_r^+= \mathcal{U}_r \cap \mathcal{U}^+$, 
$\mathcal{U}_r^0= \mathcal{U}_r \cap \mathcal{U}^0$,  and 
$\mathcal{U}_r^-= \mathcal{U}_r \cap \mathcal{U}^-$. 
Let $\mathcal{U}^{\geq 0}$ (resp. $\mathcal{U}^{\leq 0}$) be 
the $k$-subalgebra of $\mathcal{U}$ generated by $\mathcal{U}^+$ and $\mathcal{U}^0$ 
(resp. $\mathcal{U}^-$ and $\mathcal{U}^0$), and set 
$\mathcal{U}_r^{\geq 0}= \mathcal{U}_r \cap \mathcal{U}^{\geq 0}$ and 
$\mathcal{U}_r^{\leq 0}= \mathcal{U}_r \cap \mathcal{U}^{\leq 0}$.

For a  $\mathcal{U}_r^0$-module $M$, we say that a nonzero 
element $v \in M$ is a $\mathcal{U}_r^0$-weight vector of 
$\mathcal{U}_r^0$-weight $\lambda \in \mathbb{Z}$ 
if ${H \choose n} v = {\lambda \choose n} v$ for all integers $n$ with $0 \leq n \leq p^r-1$. 
Note that each $\mathcal{U}_r^0$-weight is  determined uniquely modulo $p^r$.

Let   ${\rm Fr}': \mathcal{U} \rightarrow \mathcal{U}$ be the $k$-linear map defined 
by 
$$Y^{(m)} {H \choose n} X^{(m')} \mapsto 
Y^{(pm)} {H \choose pn} X^{(pm')}$$
for $m,n,m' \in \mathbb{Z}_{\geq 0}$. 
This map is not a homomorphism of $k$-algebras, whereas its restriction to 
$\mathcal{U}^{\geq 0}$ or $\mathcal{U}^{\leq 0}$ is (for details, see  
\cite[\S1]{gros-kaneda11}). 

The following proposition plays an important role in the next section 
(for details, see \cite[Proposition 2.3]{yoshii17}). \\

\begin{Prop}\label{Multiplication}
The multiplication map 
$\mathcal{U}_1 \otimes_k {\rm Fr}'(\mathcal{U}) \rightarrow \mathcal{U}$ 
is a $k$-linear isomorphism. 
\end{Prop}
\ 

Let $\mathcal{A}$ be the $k$-subalgebra of $\mathcal{U}$ 
which is generated by $\mathcal{U}^0$ and all $Y^{(p^i)}X^{(p^i)}$ with 
 $i \geq 0$. This subalgebra is commutative and 
has  $Y^{(m)} {H \choose n} X^{(m)}$
with $m,  n \in \mathbb{Z}_{\geq 0}$ as a $k$-basis (hence consists of all elements 
of degree $0$). So  we easily see that  $\mathcal{T}_1=\mathcal{T}_2$ on $\mathcal{A}$.  For a positive integer 
$r  \in \mathbb{Z}_{>0}$, set $\mathcal{A}_r = \mathcal{A} \cap \mathcal{U}_r$. 
This subalgebra 
is generated by $\mathcal{U}_r^0$ and all $Y^{(p^i)}X^{(p^i)}$ with 
 $0 \leq i \leq r-1$, and  has 
$Y^{(m)} {H \choose n} X^{(m)}$
with $0 \leq m, n \leq p^r-1$ as a $k$-basis. 
For details on $\mathcal{A}$, see 
\cite[\S 2]{yoshii17}. 

The following proposition  is a generalization of Proposition 2.5 in \cite{yoshii17}. \\

\begin{Prop}\label{CommProp}
Let $r \in \mathbb{Z}_{>0}$ and $n \in \mathbb{Z}_{\geq 0}$. Then the elements 
$X^{(p^r n)}$ and $Y^{(p^r n)}$ in $\mathcal{U}$ commute with all elements in $\mathcal{A}_r$. 
\end{Prop}

\noindent {\itshape Proof.} We may assume that $n>0$. 
It is enough to check that all the elements 
$Y^{(s)} X^{(s)}$ and ${H \choose s}$ in $\mathcal{U}$ with   
$1 \leq s \leq p^r-1$ commute with $X^{(p^r n)}$ and $Y^{(p^r n)}$. Consider 
two equalities 
$$X^{(p^r n)} Y^{(s)} X^{(s)} = 
\sum_{i=0}^{s} Y^{(s-i)} {H-s-p^r n+2i \choose i} X^{(p^r n-i)}X^{(s)}$$
and
$$Y^{(s)} X^{(s)}Y^{(p^r n)}  = 
\sum_{i=0}^{s} Y^{(s)} Y^{(p^r n-i)}{H-s-p^r n+2i \choose i} X^{(s-i)}$$ 
in $\mathcal{U}$. In the right-hand sides of these equalities, if $i \neq 0$, we have 
$$X^{(p^r n-i)} X^{(s)}= X^{(p^r(n-1))} X^{(p^r-i)} X^{(s)} = 
{p^r-i+s \choose s} X^{(p^r(n-1))} X^{(p^r-i+s)}=0$$
and 
$$Y^{(s)}Y^{(p^r n-i)} = Y^{(s)}Y^{(p^r-i)}  Y^{(p^r(n-1))} = 
{p^r-i+s \choose s} Y^{(p^r-i+s)}Y^{(p^r(n-1))} =0,$$ 
since ${p^r-i+s \choose s}={s-i \choose s} =0$ in $\mathbb{F}_p$. Thus we have 
$X^{(p^r n)} Y^{(s)} X^{(s)} = Y^{(s)} X^{(s)}X^{(p^r n)} $ and  
$Y^{(s)} X^{(s)}Y^{(p^r n)}  = Y^{(p^r n)} Y^{(s)} X^{(s)}$. On the other hand, 
since ${2p^r n \choose l}=0$ in $\mathbb{F}_p$ for $1 \leq l \leq s$, we have  
$${H \choose s} X^{(p^r n )}= X^{(p^r n )}{H +2p^r n \choose s} 
= X^{(p^r n )} \sum_{l=0}^{s}{2p^r n \choose l}{H  \choose s-l} 
=X^{(p^r n )}{H \choose s}$$
and 
$$Y^{(p^r n )}{H \choose s} = {H +2p^r n \choose s} Y^{(p^r n )}
=  \sum_{l=0}^{s}{2p^r n \choose l}{H  \choose s-l} Y^{(p^r n )}
={H \choose s}Y^{(p^r n )}.$$
Therefore, the proposition 
follows. $\square$

\section{The elements $B^{(\bm \varepsilon)}({\bm a}, {\bm j})$}
From now until the end of this article, 
$r$ denotes a fixed positive integer unless otherwise stated. 

 In this section we shall introduce the elements $B^{(\bm \varepsilon)}({\bm a}, {\bm j})$ 
in $\mathcal{U}_r$, which include pairwise orthogonal  primitive idempotents 
whose sum is the unity $1$. For details, see also \cite[\S 4 and 5]{yoshii17},  
\cite[\S 3]{yoshii18}, and \cite[\S 3 and 4]{yoshii22}. 
 
For $a \in \mathbb{Z}$, set 
$$\mu^{(r)}_a = {H-a-1 \choose p^r-1} \in \mathcal{U}_r^0.$$
This is a   $\mathcal{U}_r^0$-weight vector of 
$\mathcal{U}_r^{0}$-weight $a$ in the $\mathcal{U}_r^0$-module  
$\mathcal{U}_r^0$. 
We have  $\mu^{(r)}_a = \mu^{(r)}_b$ if and only if $a \equiv b\ ({\rm mod}\ p^r)$, 
and all $\mu^{(r)}_a$ with $0 \leq a \leq p^r-1$ form a 
$k$-basis of $\mathcal{U}_r^0$ and 
are pairwise orthogonal 
primitive idempotents in $\mathcal{U}_r^0$ whose sum is $1$. 
It is easy to see that 
$$\mathcal{T}_1 \left( \mu_a^{(r)}\right)=\mathcal{T}_2 \left( \mu_a^{(r)}\right)
=\mu_{-a}^{(r)}$$ for any $a \in \mathbb{Z}$. We shall write 
$\mu_a^{(1)}$ simply as $\mu_a$. \\ 

\begin{Prop}\label{FormulasOfMu}
For $a \in \mathbb{Z}$, the following hold. \\ 

\noindent {\rm (i)} For any $n \in \mathbb{Z}_{\geq 0}$, we have 
$$\mu^{(r)}_a X^{(n)}= X^{(n)}\mu^{(r)}_{a-2n}\ \ \ \mbox{and}\ \ \ 
\mu^{(r)}_a Y^{(n)}= Y^{(n)}\mu^{(r)}_{a+2n}.$$
\ 

\noindent {\rm (ii)} Suppose that $r \geq 2$. 
If $a=a'+p^i a''$ with $1 \leq i \leq r-1$, $0 \leq a' \leq p^i-1$, 
and $a'' \in \mathbb{Z}$, 
we have 
$$\mu^{(r)}_a= \mu_{a'}^{(i)} {\rm Fr}'^i \left( \mu_{a''}^{(r-i)}\right).$$
\end{Prop}
\ 

\noindent For details, see \cite[\S 4]{gros-kaneda15}.
\\

Now we define a set $\mathcal{P}_{\mathbb{Z}}$ as 
$$\mathcal{P}_{\mathbb{Z}}=\mathbb{Z} \times \left\{0,1, \dots, \dfrac{p-1}{2} \right\}$$
if $p$ is odd and 
$$\mathcal{P}_{\mathbb{Z}}=\left\{ \left. \left(2i,\dfrac{1}{2}\right), (1+2i,0), (1+2i,1) 
\ \ \right|\   
i \in \mathbb{Z} \right\} \subset 
\mathbb{Z} \times \mathbb{Q}$$
if $p=2$. We also consider the subset
$$\mathcal{P} = \{(a,j) \in \mathcal{P}_{\mathbb{Z}}\ |\ 0 \leq a \leq p-1\}.$$
Clearly we have 
$$\mathcal{P} =\{ 0,1, \dots, p-1\} \times \left\{ 0,1,\dots, \dfrac{p-1}{2}\right\}$$
if $p$ is odd and $\mathcal{P} =\{ (0,1/2), (1,0),(1,1)\}$ if $p=2$. 

For an integer $n \in \mathbb{Z}$, we denote by $n\ {\rm {\bf mod}}\ p$  
a unique integer $\widehat{n}$ with 
$\widehat{n} \equiv n\ ({\rm mod}\ p)$ and $0 \leq \widehat{n} \leq p-1$. 

We classify  pairs $(a,j) \in \mathcal{P}_{\mathbb{Z}}$  
under the following four conditions: \\ \\
(A) $\widehat{a}$ is even and $(p-\widehat{a}+1)/2 \leq j \leq (p-1)/2$,\\
(B) $\widehat{a}$ is even and $0 \leq j \leq (p-\widehat{a}-1)/2$,\\
(C) $\widehat{a}$ is odd and $0 \leq j \leq (\widehat{a}-1)/2$,\\
(D) $\widehat{a}$ is odd and $(\widehat{a}+1)/2 \leq j \leq (p-1)/2$,\\ \\
where $\widehat{a}= a\ {\rm {\bf mod}}\ p$. 
Note that if $p=2$, the pairs $(2i,1/2)$, $(1+2i,0)$, and $(1+2i,1)$ in 
$\mathcal{P}_{\mathbb{Z}}$  
for $i \in \mathbb{Z}$ satisfy (B), (C), and (D) 
respectively.  Apart from them,  we also consider the following condition for 
$(a,j) \in \mathcal{P}_{\mathbb{Z}}$: \\ \\
(E) $j=0$ if $p$ is odd or $a \equiv 1\ ({\rm mod}\ 2)$ if $p=2$. \\

\begin{Def}\label{Def1}
Let $\varepsilon \in \mathbb{F}_2$ and  $(a,j) \in \mathcal{P}_{\mathbb{Z}}$, and 
set $\widehat{a}= a\ {\rm {\bf mod}}\ p$. Then 
define  nonnegative integers 
$n^{(\varepsilon)}(a,j)$ and 
$\widetilde{n}^{(\varepsilon)}(a,j)$ every condition of $(a,j)$ 
from {\rm (A)} to {\rm (D)} as follows: 
$$\begin{array}{|c||c|c|c|c|} \hline 
& & & & \\[-3mm]
(a,j) & n^{(0)}(a,j) & n^{(1)}(a,j) & \widetilde{n}^{(0)}(a,j) & \widetilde{n}^{(1)}(a,j) \\[1mm] \hline \hline 
& & & & \\[-3mm] 
{\rm (A)} & \dfrac{p-\widehat{a}-1}{2}+j   & \dfrac{3p-\widehat{a}-1}{2}-j & 
\dfrac{-p+\widehat{a}-1}{2}+j & \dfrac{p+\widehat{a}-1}{2}-j 
\\[3mm] \hline 
& & & & \\[-3mm]
{\rm (B)} & \dfrac{p-\widehat{a}-1}{2}-j & \dfrac{p-\widehat{a}-1}{2}+j & 
\dfrac{p+\widehat{a}-1}{2}-j & \dfrac{p+\widehat{a}-1}{2}+j 
\\[3mm] \hline 
& & & & \\[-3mm]
{\rm (C)} & \dfrac{2p-\widehat{a}-1}{2}-j & \dfrac{2p-\widehat{a}-1}{2}+j & 
\dfrac{\widehat{a}-1}{2}-j & \dfrac{\widehat{a}-1}{2}+j 
\\[3mm] \hline 
& & & & \\[-3mm]
{\rm (D)} & j-\dfrac{\widehat{a}+1}{2} & \dfrac{2p-\widehat{a}-1}{2}-j & 
\dfrac{\widehat{a}-1}{2}+j & \dfrac{2p+\widehat{a}-1}{2}-j 
\\[3mm] \hline 
\end{array}$$ 
\end{Def}
\ \\

\noindent {\bf Remark.} For $(a,j) \in \mathcal{P}_{\mathbb{Z}}$ and 
$\varepsilon \in \mathbb{F}_2$, we easily see the following. \\ 

\noindent {\rm (a)} $\widetilde{n}^{(\varepsilon)}(a,j)=n^{(\varepsilon)}(-a,j)$. \\ 

\noindent {\rm (b)} $0 \leq n^{(0)}(a,j) \leq n^{(1)}(a,j) \leq p-1$ and 
$$n^{(0)}(a,j)=n^{(1)}(a,j) \Longleftrightarrow 
\mbox{$(a,j)$ satisfies (E)}. $$

\noindent {\rm (c)}
$n^{(0)}(a,j)+ \widetilde{n}^{(1)}(a,j)=n^{(1)}(a,j)+ \widetilde{n}^{(0)}(a,j)=p-1$. \\

\noindent {\rm (d)} If $(a,j)$ satisfies {\rm (A)}, {\rm (C)}, or {\rm (D)}, then 
$(-a,j)$ satisfies {\rm (D)}, {\rm (B)}, or {\rm (A)} respectively. If $(a,j)$ satisfies {\rm (B)}, 
then  $(-a,j)$ satisfies {\rm (C)} if $a \not\equiv 0\ ({\rm mod}\ p)$ and {\rm (B)} if 
$a \equiv 0\ ({\rm mod}\ p)$. 

\ 

Now we recall the construction of the elements 
$B^{(\varepsilon)}(a,j) \in \mathcal{A}_1$ 
for $(a,j) \in \mathcal{P}_{\mathbb{Z}}$ and $\varepsilon \in \mathbb{F}_2$ defined in 
\cite[\S 3]{yoshii18}. 

Suppose for a moment  that $p$ is odd. 
For $\varepsilon \in \mathbb{F}_2$ and $0 \leq j \leq (p-1)/2$,  
 we define  polynomials $\psi(x), \psi_j^{(\varepsilon)}(x) \in \mathbb{F}_p [x]$  as  follows: 
$$\psi(x)=\prod_{i \in \mathbb{F}_p} \left(x- i^2\right),$$
$$\psi_0^{(0)}(x)= \psi_0^{(1)}(x)= \prod_{i \in \mathbb{F}_p^{\times}} \left(x-i^2\right),$$
$$\psi_{s}^{(0)}(x) = 2x\left(x+s^2\right) \prod_{i \in \mathbb{F}_p^{\times} 
\backslash \{ s, p-s \} } \left(x-i^2\right),  
$$
$$\psi_s^{(1)}(x) =x (x-s^2)
\prod_{i \in \mathbb{F}_p^{\times} 
\backslash \{ s, p-s \} } \left(x-i^2\right)$$
for $1 \leq s \leq (p-1)/2$  
($s$ and $p-s$ in the right-hand sides denote their images under the natural map 
$\mathbb{Z} \rightarrow \mathbb{F}_p$). 
Set 
$$B^{(\varepsilon)}(a, j) = \mu_a \cdot \psi_j^{(\varepsilon)} 
\left(\mu_a YX +\left(\dfrac{a+1}{2}\right)^2\right) \ 
\left(=   \psi_j^{(\varepsilon)} 
\left(\mu_a YX +\left(\dfrac{a+1}{2}\right)^2\right) \cdot \mu_a \right)$$
for $\varepsilon \in \mathbb{F}_2$ and $(a,j) \in \mathcal{P}_{\mathbb{Z}}$.  
Here the element $(a+1)/2$ in the right-hand side denotes $2^{-1} (a+1)$ in $\mathbb{F}_p$, 
where $2^{-1}$ is 
the inverse of $2$ in $\mathbb{F}_p$.  
Since $\psi_0^{(0)}(x)=\psi_0^{(1)}(x)$, clearly we have $B^{(0)}(a,0)=B^{(1)}(a,0)$ for any 
$a \in \mathbb{Z}$. 

In turn, suppose that  $p=2$. Then we  define 
$$B^{(0)}\left(2i,\dfrac{1}{2}\right)= \mu_0,\ \ \ B^{(1)}\left(2i,\dfrac{1}{2}\right)= 
\mu_0YX= \mu_0 XY, $$
$$B^{(0)}(1+2i,0)=B^{(1)}(1+2i,0)= \mu_1YX = \mu_1 XY + \mu_1,$$
$$B^{(0)}(1+2i,1)=  B^{(1)}(1+2i,1)= \mu_1YX + \mu_1=\mu_1XY$$
for any $i \in \mathbb{Z}$.
\\

Let $(a,j) \in \mathcal{P}_{\mathbb{Z}}$. 
For $i \in \mathbb{Z}$ and $n \in \{0,1, \dots, p-1 \}$, define 
$\gamma_{i}(a,j)$, $\widetilde{\gamma}_{i}(a,j)$, $\beta_{n}(a,j)$, and 
$\widetilde{\beta}_{n}(a,j)$ in $\mathbb{F}_p$ as follows: 
$$\gamma_{i}(a,j)= j^2- \left( \dfrac{a+1}{2} \right)^2 -i (i+a+1) \ \left( =
 j^2- \left( \dfrac{a+1}{2} +i \right)^2 \right),$$
$$\widetilde{\gamma}_{i}(a,j)= \gamma_{i}(-a,j),$$
$$\beta_{n}(a,j)= \prod_{i=0}^{n-1}\gamma_{i}(a,j),$$
$$\widetilde{\beta}_{n}(a,j) = \beta_{n}(-a,j) \ \left( 
=\prod_{i=0}^{n-1}\widetilde{\gamma}_{i}(a,j) \right).$$  
Here if $p=2$, $\gamma_{i}(a,j)$ is defined by regarding the right-hand side (which is 
an integer in this situation) as the image under the natural map
$\mathbb{Z} \rightarrow \mathbb{F}_2$. If $p$ is odd, 
$\gamma_{i}(a,j)$ is defined by regarding the integers $i$, $j$, and $a+1$ in the 
right-hand side as the images 
under the natural map $\mathbb{Z} \rightarrow \mathbb{F}_p$. \\ 

\begin{Prop}\label{PropertiesOfB1}
For $(a,j) \in \mathcal{P}_{\mathbb{Z}}$ and $\varepsilon \in \mathbb{F}_2$, 
the following hold. \\

\noindent {\rm (i)} $B^{(\varepsilon)}(a,j)$ has $\mathcal{U}_1^0$-weight $a$: 
$H B^{(\varepsilon)}(a,j) = a B^{(\varepsilon)}(a,j)$. \\

\noindent {\rm (ii)} The elements 
$B^{(0)}(a,j)$ with $(a,j) \in \mathcal{P}$ are pairwise orthogonal 
primitive idempotents in $\mathcal{U}_1$ whose sum is the unity $1 \in \mathcal{U}_1$. \\
\ 

\noindent {\rm (iii)} $B^{(\varepsilon)}(a,j)$ can be written as 
$$B^{(\varepsilon)}(a,j)= \mu_a \sum_{m=n^{(\varepsilon)}(a,j)}^{p-1} 
c^{(\varepsilon)}_m(a,j) Y^m X^m = 
\mu_a \sum_{m=\widetilde{n}^{(\varepsilon)}(a,j)}^{p-1} \widetilde{c}^{(\varepsilon)}_m(a,j) X^m Y^m$$
for some $c^{(\varepsilon)}_m(a,j), \widetilde{c}^{(\varepsilon)}_m(a,j) \in \mathbb{F}_p$ with 
$c^{(\varepsilon)}_{n^{(\varepsilon)}(a,j)}(a,j) \neq 0$ and 
$\widetilde{c}^{(\varepsilon)}_{\widetilde{n}^{(\varepsilon)}(a,j)}(a,j) \neq 0$. \\

\noindent {\rm (iv)} We have 
$$YXB^{(0)}(a,j)=\gamma_0(a,j) B^{(0)}(a,j) + 4j^2 B^{(1)}(a,j),$$
$$XYB^{(0)}(a,j)=\widetilde{\gamma}_0(a,j) B^{(0)}(a,j) + 4j^2 B^{(1)}(a,j),$$
$$YXB^{(1)}(a,j)=\gamma_0(a,j) B^{(1)}(a,j),$$
$$XYB^{(1)}(a,j)=\widetilde{\gamma}_0(a,j) B^{(1)}(a,j).$$
\ 
 
\noindent {\rm (v)} We have  
$$B^{(0)}(a,j)=B^{(1)}(a,j) \Longleftrightarrow  
\mbox{$(a,j)$ satisfies {\rm (E)}}.$$
\ 

\noindent {\rm (vi)}  $B^{(\varepsilon)}(a,j)B^{(0)}(a,j)=
B^{(0)}(a,j)B^{(\varepsilon)}(a,j)=B^{(\varepsilon)}(a,j)$. \\ 

\noindent {\rm (vii)} For $(a',j') \in \mathcal{P}_{\mathbb{Z}}$, 
we have  
$$B^{(\varepsilon)}(a,j)=B^{(\varepsilon)}(a',j') \Longleftrightarrow 
\mbox{$a \equiv a'\ ({\rm mod}\ p)$ and $j=j'$}.$$ 
\

\noindent {\rm (viii)} If $p$ is odd, then 
$$X B^{(\varepsilon)}(a,j) = B^{(\varepsilon)}(a+2,j) X\ \ \ \mbox{and}\ \ \ \ 
Y B^{(\varepsilon)}(a,j) = B^{(\varepsilon)}(a-2,j) Y. $$
If $p=2$, then 
$$
XB^{(\varepsilon)}\left(0, \dfrac{1}{2}\right) =
B^{(\varepsilon)}\left(0, \dfrac{1}{2}\right) X, \ 
YB^{(\varepsilon)}\left(0, \dfrac{1}{2}\right)=
B^{(\varepsilon)}\left(0, \dfrac{1}{2}\right) Y,
$$
$$XB^{(\varepsilon)}(1,0)= B^{(\varepsilon)}(1,1)X, \ 
YB^{(\varepsilon)}(1,0)= B^{(\varepsilon)}(1,1)Y,$$
$$
XB^{(\varepsilon)}(1,1)= B^{(\varepsilon)}(1,0)X, \     
YB^{(\varepsilon)}(1,1)= B^{(\varepsilon)}(1,0)Y.$$ 
\end{Prop}

\noindent {\itshape Proof.} (i) is clear by the definition of $B^{(\varepsilon)}(a,j)$. 
For (ii), see \cite[Proposition 4.5]{yoshii17} (see also \cite[Theorem 1]{seligman03}). 
For (iii), see \cite[Lemma 3.3]{yoshii18}. For (viii), see \cite[Proposition 4.2]{yoshii22}. 
(iv) follows from direct calculation using the definition of $B^{(\varepsilon)}(a,j)$. 
(v) follows from the definition of 
$B^{(\varepsilon)}(a,j)$, (b) in the remark of Definition \ref{Def1}, and (iii).  
We shall show (vi). For $\varepsilon=0$, it follows from (ii). So assume that 
$\varepsilon=1$. Since both $B^{(0)}(a,j)$ and $B^{(1)}(a,j)$ lie in $\mathcal{A}_1$, 
they are commutative. Note that $\psi_j^{(0)}(j^2)=1$ in $\mathbb{F}_p$ (see the proof of 
\cite[Lemma 4.4]{yoshii17}). By (iv), we obtain 
\begin{align*}
B^{(0)}(a, j) B^{(1)}(a,j)&= \mu_a \cdot \psi_j^{(0)} 
\left(\mu_a YX +\left(\dfrac{a+1}{2}\right)^2\right) B^{(1)}(a,j)\\
&= \psi_j^{(0)} 
\left(\gamma_0(a,j) +\left(\dfrac{a+1}{2}\right)^2\right) B^{(1)}(a,j) \\
&=\psi_j^{(0)} 
(j^2) B^{(1)}(a,j) \\
&=B^{(1)}(a,j)
\end{align*}
and (vi) follows. Finally, as for (vii), the `if' part is clear and the `only if' part 
easily follows from (ii) and (vi). 
$\square$
\\

For $(a,j) \in \mathcal{P}_{\mathbb{Z}}$, if $(a,j)$ satisfies 
{\rm (A)} or {\rm (C)}, then define an integer $s(a,j)$ as 
$$s(a,j)=\dfrac{p-(a \ {\rm {\bf mod}}\ p)+1}{2}$$
 if $p$ is odd and 
$a\ {\rm {\bf mod}}\ p$ is even,  
$$s(a,j)=\dfrac{p-(a\ {\rm {\bf mod}}\ p)}{2}$$ if both $p$ and $a\ {\rm {\bf mod}}\ p$ are  odd, 
and  $s(a,j)=1$ if $p=2$. 

For $\varepsilon \in \mathbb{F}_2$ and $(a,j) \in \mathcal{P}_{\mathbb{Z}}$, if  we write 
$$B^{(\varepsilon)}(a,j)= \mu_a \sum_{m=n^{(\varepsilon)}(a,j)}^{p-1} 
c^{(\varepsilon)}_m(a,j) Y^m X^m$$ 
following Proposition \ref{PropertiesOfB1} (iii), then    
define $Z^{(\varepsilon)} \left(z; (a,j) \right) \in \mathcal{U}$ for $z \in \mathcal{U}$ as 
$$Z^{(\varepsilon)} \left(z; (a,j) \right)
= \mu_a \sum_{m=n^{(\varepsilon)}(a,j)}^{p-1} 
c^{(\varepsilon)}_m(a,j) Y^m X^{m-s(a,j)} {\rm Fr'}(z)X^{s(a,j)}$$
if $(a,j)$ satisfies (A) or (C), and
$$Z^{(\varepsilon)} \left(z; (a,j) \right)
={\rm Fr'}(z) B^{(\varepsilon)}(a,j)
\left( =B^{(\varepsilon)}(a,j){\rm Fr'}(z) \right)$$
if $(a,j)$ satisfies (B) or (D).  \\

\begin{Prop}\label{PropertiesOfZ} 
Let $(a,j) \in \mathcal{P}_{\mathbb{Z}}$ and $\varepsilon \in \mathbb{F}_2$. 
The following hold. \\ 

\noindent {\rm (i)} The map $Z^{(\varepsilon)} \left(-; (a,j) \right) : \mathcal{U} \rightarrow \mathcal{U}, \ 
z \mapsto Z^{(\varepsilon)} \left(z; (a,j) \right)$ is  $k$-linear and injective. \\

\noindent {\rm (ii)} For  $z \in \mathcal{U}$,  there is  
an element $z' \in \mathcal{U}$ which is independent of 
$\varepsilon$ such that $$Z^{(\varepsilon)} \left(z; (a,j) \right)
={\rm Fr'}(z') B^{(\varepsilon)}(a,j)=B^{(\varepsilon)}(a,j){\rm Fr'}(z').$$
Then we also have 
$$z=0 \Longleftrightarrow {\rm Fr}'(z')=0 \Longleftrightarrow z'=0$$
and 
$$z \in \mathcal{A} \Longleftrightarrow Z^{(\varepsilon)} \left(z; (a,j) \right) \in \mathcal{A}
\Longleftrightarrow z' \in \mathcal{A}.$$
\ 

\noindent {\rm (iii)} Let $u$ be an element of the $k$-subalgebra of $\mathcal{U}$ generated by all $X^{(p^i)}$ and $Y^{(p^i)}$ 
with $i \in \mathbb{Z}_{> 0}$ (and the unity $1 \in \mathcal{U}$). 
Then we have 
$u Z^{(\varepsilon)} \left(z; (a,j) \right)= 
Z^{(\varepsilon)} \left( {\rm Fr}(u)z; (a,j) \right)$. \\ 

\noindent {\rm (iv)} We have 
$$Z^{(\varepsilon)} \left(z_1; (a,j) \right) Z^{(0)} \left(z_2; (a,j) \right) 
=Z^{(0)} \left(z_1; (a,j) \right) Z^{(\varepsilon)} \left(z_2; (a,j) \right) 
= Z^{(\varepsilon)} \left(z_1z_2; (a,j) \right)$$ for any $z_1,z_2 \in \mathcal{U}$. \\ 

\noindent {\rm (v)} 
For a nonzero element 
$z \in \mathcal{U}$, we have 
$$Z^{(0)} \left(z; (a,j) \right)=Z^{(1)} \left(z; (a,j) \right) \Longleftrightarrow  
\mbox{$(a,j)$ satisfies {\rm (E)}}.$$ 
\ 

\noindent {\rm (vi)} 
For $(a',j') \in \mathcal{P}_{\mathbb{Z}}$ and a nonzero element 
$z \in \mathcal{U}$, we have  
$$Z^{(\varepsilon)} \left(z; (a,j) \right)=Z^{(\varepsilon)} 
\left(z; (a',j') \right) \Longleftrightarrow 
\mbox{$a \equiv a'\ ({\rm mod}\ p)$ and $j=j'$}.$$ 
\end{Prop}

\noindent {\itshape Proof.} Since the map ${\rm Fr}'$ is $k$-linear, the linearity in (i) is 
clear. For $n_1, n_2, n_3 \in \mathbb{Z}_{\geq 0}$, we have 
$$
Z^{(\varepsilon)} \left( Y^{(n_1)} {H \choose n_2} X^{(n_3)} ; (a,j) \right)  
 = B^{(\varepsilon)}(a,j) {\rm Fr}' \left( Y^{(n_1)} \left( {H \choose n_2 } + 
{H \choose n_2-1 } \right) X^{(n_3)} \right)   
$$
if $(a,j)$ satisfies (A) or (C) (for convenience, we define 
${H \choose t}=0$ if $t<0$). This equality is described in the proof of 
\cite[Lemma 5.3]{yoshii17} when $\varepsilon=0$, but 
it also holds even if  $\varepsilon=1$.  This fact together with Proposition \ref{Multiplication} 
and the injectivity of ${\rm Fr}'$ easily implies the injectivity of 
$Z^{(\varepsilon)} \left(-; (a,j) \right)$ and (ii). For (iii), see \cite[Lemma 3.5]{yoshii18} 
and the paragraph just after (see also the proof of \cite[Proposition 5.4 (iv)]{yoshii17}). 
(iv) is obtained by multiplying both sides of the equality 
$$Z^{(0)} \left(z_1; (a,j) \right) Z^{(0)} \left(z_2; (a,j) \right) 
= Z^{(0)} \left(z_1z_2; (a,j) \right)$$
(see \cite[Proposition 5.4 (ii)]{yoshii17}) by $B^{(\varepsilon)}(a,j)$. (v) follows from 
(ii), Proposition \ref{PropertiesOfB1} (v), and Proposition \ref{Multiplication}. Finally, as 
for (vi), the `if' part is clear and the `only if' part follows from (ii), Proposition 
\ref{Multiplication}, and Proposition \ref{PropertiesOfB1} (vii). $\square$
\  

Consider an $r$-tuple 
$\left((a_i, j_i)\right)_{i=0}^{r-1} =\left( (a_0,j_0),\dots, (a_{r-1},j_{r-1}) \right)
\in \mathcal{P}_{\mathbb{Z}}^r$ 
of  pairs 
$(a_i,j_i) \in \mathcal{P}_{\mathbb{Z}}$ $(0 \leq i \leq r-1)$. For convenience 
we shall write it as 
$$\left((a_0,  \dots, a_{r-1}),(j_0,  \dots, j_{r-1})\right)$$
or $({\bm a},{\bm j})$ with ${\bm a}=(a_0,  \dots, a_{r-1})$ and 
${\bm j}=(j_0,  \dots, j_{r-1})$.

For ${\bm \varepsilon}=(\varepsilon_0, \dots, \varepsilon_{r-1}) \in \mathbb{F}_2^r$,  
$({\bm a},{\bm j})= \left((a_i, j_i)\right)_{i=0}^{r-1} \in \mathcal{P}_{\mathbb{Z}}^r$ and 
$z \in \mathcal{U}$,  
we define an element 
$Z^{({\bm \varepsilon})} \left( z ; ({\bm a}, {\bm j}) \right)
\in \mathcal{U}$ inductively 
as $Z^{({\bm \varepsilon})} \left( z ; ({\bm a}, {\bm j}) \right) = 
Z^{(\varepsilon_0)}\left( z; (a_0,j_0) \right)$ if $r=1$, and 
$$Z^{({\bm \varepsilon})} \left( z ; ({\bm a}, {\bm j}) \right)=
Z^{(\varepsilon_0)} \left( 
Z^{({\bm \varepsilon}')}
\left(z; \left( {\bm a}', {\bm j}' \right) \right); (a_0,j_0) \right) 
$$
if $r \geq 2$, where ${\bm \varepsilon}'=(\varepsilon_1, \dots, \varepsilon_{r-1})$ and 
$\left({\bm a}', {\bm j}'\right)=((a_i,j_i))_{i=1}^{r-1}$. Then set 
$B^{({\bm \varepsilon})}({\bm a}, {\bm j})=Z^{({\bm \varepsilon})} 
\left( 1 ; ({\bm a}, {\bm j}) \right)
\in \mathcal{A}_r$. \\

In this paper, we need a lot of notation to prove a main result. So we collect them here 
(see also \cite[Definition 4.1]{yoshii22} and its remark). \\

\begin{Def}\label{Def2}
Let $(a,j) \in \mathcal{P}_{\mathbb{Z}}$ and 
$({\bm a}, {\bm j})= \left( (a_i,j_i)\right)_{i=0}^{r-1} \in \mathcal{P}_{\mathbb{Z}}^r$. \\

\noindent {\rm (1)} Let $\iota : \mathcal{P}_{\mathbb{Z}} \rightarrow \mathbb{Z}$ 
be a function  defined as 
$$\iota(a,j) = 
\left\{ \begin{array}{ll} 
a\ {\rm {\bf mod}}\ p-p & {\mbox{if $(a,j)$ 
satisfies {\rm (A)} or {\rm (C)},}} \\
a\ {\rm {\bf mod}}\ p   & {\mbox{if $(a,j)$ 
satisfies {\rm (B)} or {\rm (D)}}}
\end{array} \right. .$$
In this paper, set $b=\iota(a,j)$ and $b_i=\iota(a_i,j_i)$ for 
$0 \leq i \leq r-1$ and we use the notation 
unless otherwise stated. \\

\noindent {\rm (2)} We regard $\mathbb{F}_2^r$ and 
$\mathbb{Z}^r$ as  additive groups induced by 
the additions in $\mathbb{F}_2$ and $\mathbb{Z}$ respectively. 
Let ${\bm e}_i$ denote an element of 
$\mathbb{F}_2^r$ or $\mathbb{Z}^r$with $1$ in the $i$-th entry and $0$ elsewhere. 
Define two elements ${\bm 0}_r$ and ${\bm 1}_r$ in 
$\mathbb{F}_2^r$ or $\mathbb{Z}^r$ as 
$${\bm 0}_r= (0,\dots, 0),$$
$${\bm 1}_r= (1,\dots, 1).$$
We usually write them as ${\bm 0}$ and ${\bm 1}$ respectively unless confusion occurs.  
Moreover, for ${\bm \varepsilon}=(\varepsilon_0, \dots, \varepsilon_{r-1})$, 
$\widetilde{\bm \varepsilon}=(\widetilde{\varepsilon}_0, \dots, \widetilde{\varepsilon}_{r-1})
 \in \mathbb{F}_2^r$,  
define   
${\bm \varepsilon} \leq \widetilde{\bm \varepsilon}$ 
 if $\varepsilon_i \leq \widetilde{\varepsilon}_i$ for each $i$,   
regarding $\varepsilon_i$ and $\widetilde{\varepsilon}_i$ as the corresponding integers 
(i.e. $0$ or $1$ in $\mathbb{Z}$). This 
gives a partial order in $\mathbb{F}_2^{r}$. \\ 

\noindent {\rm (3)} A subset $\mathcal{X}_r({\bm a},{\bm j})$ of 
$\mathbb{F}_2^r$ is defined as 
$$\mathcal{X}_r({\bm a}, {\bm j}) = 
\{ (\varepsilon_0, \dots, \varepsilon_{r-1})\in \mathbb{F}_2^r\ |\ 
\varepsilon_i=0\ \mbox{whenever $(a_i,j_i)$ satisfies {\rm (E)}}\}.$$
\

\noindent {\rm (4)} For ${\bm \varepsilon} \in \mathcal{X}_r({\bm a},{\bm j})$, define a 
 subset  $\Theta_r\left( ({\bm a}, {\bm j}), {\bm \varepsilon}\right)$  of 
$\mathbb{F}_2^r \times \mathbb{Z}^r$ as 
$$\Theta_r\left( ({\bm a}, {\bm j}), {\bm \varepsilon}\right)=
\left\{ \left( {\bm \theta}, {\bm t}({\bm \theta})\right)\ \left|\  
\begin{array}{l}
\mbox{${\bm \varepsilon} \leq {\bm \theta} \in \mathcal{X}_r({\bm a},{\bm j})$ and} \\
\mbox{$-\widetilde{n}^{(\theta_i+1)}(a_i,j_i) \leq t_i(\theta_i) \leq n^{(\theta_i+1)}(a_i,j_i)$ 
for each $i$}
\end{array}
\right. \right\},$$
where ${\bm \theta}= (\theta_0, \dots, \theta_{r-1}) \in \mathbb{F}_2^r$ and 
${\bm t}({\bm \theta})= \left( t_0(\theta_0), \dots, t_{r-1}(\theta_{r-1})\right) 
\in \mathbb{Z}^r$.  
From now on we adopt such notation for the entries of 
${\bm \theta}$ and ${\bm t}({\bm \theta})$ with respect to an element 
$\left( {\bm \theta}, {\bm t}({\bm \theta})\right)$ in 
$\Theta_r\left( ({\bm a}, {\bm j}), {\bm \varepsilon}\right)$ 
unless otherwise stated. 
\\ 

\noindent {\rm (5)} For $i \in \mathbb{Z}_{\geq 0}$ and 
$t \in \mathbb{Z}$, define an element 
$u^{(i, t)}$ in $\mathcal{U}$ as 
$$u^{(i, t)}=\left\{ \begin{array}{ll} 
 {X^{(p^i)t}} & {\mbox{if $t \geq 0$,}} \\ 
 {\left(Y^{(p^i)}\right)^{-t}} & {\mbox{if $t < 0$}}
\end{array} \right. .$$
Moreover, for   
${\bm \varepsilon} =(\varepsilon_0, \dots, \varepsilon_{r-1}) 
\in \mathbb{F}_2^r$ and ${\bm t}= (t_0 ,\dots, t_{r-1}) \in \mathbb{Z}^r$, 
define an element 
$B^{({\bm \varepsilon})}\left( ({\bm a}, {\bm j}); {\bm t} \right)$ in $\mathcal{U}_r$ as 
$$B^{({\bm \varepsilon})}\left( ({\bm a}, {\bm j}); {\bm t} \right)
=u^{(0,t_0)} u^{(1, t_1)} \cdots u^{({r-1}, t_{r-1})}
B^{({\bm \varepsilon})}({\bm a}, {\bm j}).$$

\noindent {\rm (6)} For ${\bm \varepsilon} \in \mathcal{X}_r({\bm a},{\bm j})$, 
define a subset  
$\mathcal{B}_r\left( ({\bm a}, {\bm j}), {\bm \varepsilon}\right)$  of 
$\mathcal{U}_r$ as 
$$\mathcal{B}_r\left( ({\bm a}, {\bm j}), {\bm \varepsilon}\right)
= \left\{ \left. B^{({\bm \theta})}\left(({\bm a}, {\bm j}); {\bm t}({\bm \theta})\right)
\ \right|\ \left({\bm \theta}, {\bm t}({\bm \theta})\right) \in 
\Theta_{r}\left( ({\bm a}, {\bm j}), {\bm \varepsilon}\right)\right\}.$$
\end{Def} 
\ 

\noindent {\bf Remark.} We have 
\begin{align*}
b+2n^{(0)}(a,j) &= 2 \widetilde{n}^{(0)}(a,j)-b \\
&= \left\{
\begin{array}{ll}
2j-1 & \mbox{if $(a,j)$ satisfies (A) or (D)}, \\
p-2j-1 & \mbox{if $(a,j)$ satisfies (B) or (C)}
\end{array}
\right.,
\end{align*}
where $b$ is the integer defined in (1). \\

The following proposition is a generalization of Proposition \ref{PropertiesOfB1} (iv). \\

\begin{Prop}\label{FormulasOfB1}
For $\varepsilon \in \mathbb{F}_2$ and $(a,j) \in \mathcal{P}_{\mathbb{Z}}$,  
we have 
$$Y^sX^sB^{(0)}(a,j) = \beta_s(a,j)B^{(0)}(a,j)+4j^2\sum_{i=0}^{s-1}
\dfrac{\beta_{s}(a,j)}{\gamma_i(a,j)} B^{(1)}(a,j), $$ 
$$X^tY^tB^{(0)}(a,j) = \widetilde{\beta}_t(a,j)B^{(0)}(a,j)+4j^2\sum_{i=0}^{t-1}
\dfrac{\widetilde{\beta}_{t}(a,j)}{\widetilde{\gamma}_i(a,j)} B^{(1)}(a,j)$$
if $0 \leq s \leq n^{(0)}(a,j)$ and $0 \leq t \leq \widetilde{n}^{(0)}(a,j)$,  
$$Y^s X^s B^{(0)}(a,j) = 4j^2 
\left( \prod_{i=0,\ i \neq {n}^{(0)}(a,j)}^{s-1}\gamma_i(a,j) \right) B^{(1)}(a,j),$$
$$X^t Y^t B^{(0)}(a,j) = 4j^2 
\left( \prod_{i=0,\ i \neq \widetilde{n}^{(0)}(a,j)}^{t-1}
\widetilde{\gamma}_i(a,j) \right) B^{(1)}(a,j)$$
if $n^{(0)}(a,j) < s \leq p-1$ and $\widetilde{n}^{(0)}(a,j) < t \leq p-1$, and 
$$Y^sX^sB^{(1)}(a,j) = \beta_s(a,j)B^{(1)}(a,j),$$
$$X^tY^tB^{(1)}(a,j) = \widetilde{\beta}_t(a,j)B^{(1)}(a,j)$$ 
if $0 \leq s,t \leq p-1$. 
\end{Prop}

\noindent {\itshape Proof.} See \cite[Proposition 4.3]{yoshii22}. $\square$ \\

Now we describe some properties of the elements 
$B^{(\bm{\varepsilon} )}({\bm a},{\bm j})$ and 
$Z^{({\bm \varepsilon})} \left(z; ({\bm a},{\bm j}) \right)$ for 
$({\bm a},{\bm j}) \in \mathcal{P}_{\mathbb{Z}}^r$,  
${\bm \varepsilon} \in \mathbb{F}_2^r$, and  $z \in \mathcal{U}$. \\

\begin{Prop}\label{PropertiesOfBZ} 
Let $({\bm a},{\bm j}) = \left( (a_i,j_i)\right)_{i=0}^{r-1} \in \mathcal{P}_{\mathbb{Z}}^r$ and 
${\bm \varepsilon}=(\varepsilon_0, \dots, \varepsilon_{r-1})
\in \mathbb{F}_2^r$. 
The following hold. \\ 

\noindent {\rm (i)}  $B^{(\bm{\varepsilon} )}({\bm a},{\bm j})$ has 
 $\mathcal{U}_r^0$-weight $\sum_{i=0}^{r-1}p^i b_i $ (recall that 
$b_i= \iota(a_i,j_i)$). \\ 

\noindent {\rm (ii)} For $0 \leq i \leq r-1$, we have 
$$Y^{(p^i)s}X^{(p^i)s}B^{({\bm \varepsilon})}({\bm a},{\bm j}) 
= \beta_s(a_i,j_i)B^{({\bm \varepsilon})}({\bm a},{\bm j})+4j_i^2\sum_{l=0}^{s-1}
\dfrac{\beta_{s}(a_i,j_i)}{\gamma_l(a_i,j_i)} 
B^{({\bm \varepsilon}+{\bm e}_{i+1})}({\bm a},{\bm j}), $$ 
$$X^{(p^i)t}Y^{(p^i)t}B^{({\bm \varepsilon})}({\bm a},{\bm j}) 
= \widetilde{\beta}_t(a_i,j_i)B^{({\bm \varepsilon})}({\bm a},{\bm j})+4j_i^2\sum_{l=0}^{t-1}
\dfrac{\widetilde{\beta}_{t}(a_i,j_i)}{\widetilde{\gamma}_l(a_i,j_i)} 
B^{({\bm \varepsilon}+{\bm e}_{i+1})}({\bm a},{\bm j})$$
if $\varepsilon_i=0$, $0 \leq s \leq n^{(0)}(a_i,j_i)$, and $0 \leq t \leq \widetilde{n}^{(0)}(a_i,j_i)$,  
$$Y^{(p^i)s} X^{(p^i)s} B^{({\bm \varepsilon})}({\bm a},{\bm j}) = 4j_i^2 
\left( \prod_{l=0,\ l \neq {n}^{(0)}(a_i,j_i)}^{s-1}\gamma_l(a_i,j_i) \right) 
B^{({\bm \varepsilon}+{\bm e}_{i+1})}({\bm a},{\bm j}),$$
$$X^{(p^i)t} Y^{(p^i)t} B^{({\bm \varepsilon})}({\bm a},{\bm j}) = 4j_i^2 
\left( \prod_{l=0,\ l \neq \widetilde{n}^{(0)}(a_i,j_i)}^{t-1}\widetilde{\gamma}_l(a_i,j_i) \right) 
B^{({\bm \varepsilon}+{\bm e}_{i+1})}({\bm a},{\bm j})$$
if $\varepsilon_i=0$, $n^{(0)}(a_i,j_i) < s \leq p-1$, and 
$\widetilde{n}^{(0)}(a_i,j_i) < t \leq p-1$, and 
$$Y^{(p^i)s}X^{(p^i)s}B^{({\bm \varepsilon})}({\bm a},{\bm j}) = 
\beta_s(a_i,j_i)B^{({\bm \varepsilon})}({\bm a},{\bm j}),$$
$$X^{(p^i)t}Y^{(p^i)t}B^{({\bm \varepsilon})}({\bm a},{\bm j}) = 
\widetilde{\beta}_t(a_i,j_i)B^{({\bm \varepsilon})}({\bm a},{\bm j})$$
if $\varepsilon_i=1$ and $0 \leq s,t \leq p-1$. \\ 

\noindent {\rm (iii)} The map $Z^{({\bm \varepsilon})} 
\left(-; ({\bm a},{\bm j}) \right) : \mathcal{U} \rightarrow \mathcal{U}, \ 
z \mapsto Z^{({\bm \varepsilon})} \left(z; ({\bm a},{\bm j}) \right)$ is  
$k$-linear and injective. \\ 

\noindent {\rm (iv)} Let $u$ be an element of the $k$-subalgebra of $\mathcal{U}$ generated by all $X^{(p^i)}$ and $Y^{(p^i)}$ 
with $i \geq r$ (and the unity $1 \in \mathcal{U}$). 
Then we have 
$u Z^{({\bm \varepsilon})} \left(z; ({\bm a},{\bm j}) \right)= 
Z^{({\bm \varepsilon})} \left( {\rm Fr}^r(u)z; ({\bm a},{\bm j}) \right)$ 
for any $z \in \mathcal{U}$. \\

\noindent {\rm (v)} For $r' \in \mathbb{Z}_{>0}$ and $a' \in \mathbb{Z}$, we have 
$Z^{({\bm \varepsilon})} \left( \mu_{a'}^{(r')}z; ({\bm a},{\bm j}) \right)= 
\mu_{\sum_{i=0}^{r-1} p^i b_i +p^r a'}^{(r+r')}
Z^{({\bm \varepsilon})} \left( z; ({\bm a},{\bm j}) \right)$ 
for any $z \in \mathcal{U}$. \\ 
 
\noindent {\rm (vi)} For  $z \in \mathcal{U}$,  there is  
an element $z' \in \mathcal{U}$ which is independent of 
${\bm \varepsilon}$ such that $$Z^{({\bm \varepsilon})} \left(z; ({\bm a},{\bm j}) \right)
={\rm Fr'}^r(z') B^{({\bm \varepsilon})}({\bm a},{\bm j})=
B^{({\bm \varepsilon})}({\bm a},{\bm j}){\rm Fr'}^r(z').$$
Then we also have 
$$z=0 \Longleftrightarrow {\rm Fr}'^r(z')=0 \Longleftrightarrow z'=0$$
and 
$$z \in \mathcal{A} \Longleftrightarrow Z^{({\bm \varepsilon})} 
\left(z; ({\bm a},{\bm j}) \right) \in \mathcal{A}
\Longleftrightarrow z' \in \mathcal{A}.$$  
\ 

\noindent {\rm (vii)} We have 
$$Z^{({\bm \varepsilon})} \left(z_1; ({\bm a},{\bm j}) \right) 
Z^{({\bm 0})} \left(z_2; ({\bm a},{\bm j}) \right) 
=Z^{({\bm 0})} \left(z_1; ({\bm a},{\bm j}) \right) 
Z^{({\bm \varepsilon})} \left(z_2; ({\bm a},{\bm j}) \right) 
= Z^{({\bm \varepsilon})} \left(z_1z_2; ({\bm a},{\bm j}) \right)$$ 
for any $z_1,z_2 \in \mathcal{U}$. \\ 

\noindent {\rm (viii)} 
For $\widetilde{\bm \varepsilon}=(\widetilde{\varepsilon}_0, \dots, 
\widetilde{\varepsilon}_{r-1}) \in \mathbb{F}_2^r$ and 
a nonzero element $z \in \mathcal{U}$, we have 
$$Z^{(\bm{\varepsilon} )}\left( z;({\bm a},{\bm j}) \right)= 
Z^{(\widetilde{\bm{\varepsilon}} )} \left( z;({\bm a},{\bm j}) \right) \Longleftrightarrow   
\mbox{$\varepsilon_i = \widetilde{\varepsilon}_i$ whenever 
$(a_i,j_i)$ does not satisfy {\rm (E)}}.$$ 
\ 

\noindent {\rm (ix)} For $\left(\widetilde{\bm a},\widetilde{\bm j}\right)= 
\left( \left(\widetilde{a}_i, \widetilde{j}_i \right)\right)_{i=0}^{r-1} 
\in \mathcal{P}_{\mathbb{Z}}^r$ and a nonzero $z \in \mathcal{U}$, 
we have  
$$Z^{({\bm \varepsilon})} \left( z ; ({\bm a}, {\bm j}) \right)=
Z^{({\bm \varepsilon})} \left( z ; \left(\widetilde{\bm a}, \widetilde{\bm j}\right) \right)
\Longleftrightarrow  \mbox{$a_i \equiv \widetilde{a}_i\ 
({\rm mod}\ p)$ and $j_i=\widetilde{j}_i$ for each $i$}.$$ 
\

\noindent {\rm (x)} The elements $B^{({\bm 0})}({\bm a}, {\bm j})$ with $({\bm a}, {\bm j}) \in \mathcal{P}^{r}$ 
are pairwise orthogonal primitive 
idempotents in $\mathcal{U}_r$ whose sum is the unity $1 \in \mathcal{U}_r$. 
\end{Prop} 

\noindent {\itshape Proof.} (x) is proved in \cite[Proposition 5.5 (iii)]{yoshii17}. 
(i) can be proved as in \cite[Proposition 5.5(i)]{yoshii17}.    (ii), (iii), (iv), (vii), (viii), and 
(ix) are 
easily proved by induction on $r$ using Propositions \ref{PropertiesOfZ} and 
\ref{FormulasOfB1}. For (v), since 
$$Z^{({\bm \varepsilon})} \left( \mu_{a'}^{(r')}z; ({\bm a},{\bm j}) \right)=
Z^{({\bm 0})} \left( \mu_{a'}^{(r')}; ({\bm a},{\bm j}) \right)
Z^{({\bm \varepsilon})} \left( z; ({\bm a},{\bm j}) \right),$$
it is enough to 
check that $$Z^{({\bm 0})} \left( \mu_{a'}^{(r')}; ({\bm a},{\bm j}) \right)= 
\mu_{\sum_{i=0}^{r-1} p^i b_i +p^r a'}^{(r+r')}
B^{({\bm 0})} ({\bm a},{\bm j}).$$ 
It is easy to show it by induction on $r$, since we can show that 
$$Z^{(0)} \left(  \mu_{\sum_{i=1}^{r-1}p^{i-1} b_i +p^{r-1}a'}^{(r+r'-1)}; (a_0,j_0) 
\right)
= \mu_{\sum_{i=0}^{r-1} p^i b_i +p^r a'}^{(r+r')}
B^{(0)} ({a}_0,{j}_0)$$ 
as in the proof of \cite[Proposition 5.5 (i)]{yoshii17}. 
Finally, we shall prove (vi). Note that the element $z \in \mathcal{U}$ lies in 
$\mathcal{U}_{r'}$ for some $r' \in \mathbb{Z}_{>0}$. 
Since the elements $Y^{(n_1)} \mu_{n_2}^{(r')} X^{(n_3)}$ 
with $0 \leq n_1,n_2,n_3 \leq p^{r'}-1$ form a 
$k$-basis of $\mathcal{U}_{r'}$, we may assume that $z$ is a basis element 
$Y^{(n_1)} \mu_{n_2}^{(r')} X^{(n_3)}$. Set
$$n_2'= 
\left\{ \begin{array}{ll}
n_2 & \mbox{if $\sum_{i=0}^{r-1}p^i b_i \geq 0$,} \\
n_2-1 & \mbox{if $\sum_{i=0}^{r-1}p^i b_i < 0$} 
\end{array}  
\right..$$
By Propositions \ref{PropertiesOfBZ} (iv), (v), \ref{FormulasOfMu} (ii), and  
\ref{CommProp} 
we have 
\begin{align*}
Z^{({\bm \varepsilon})} \left(Y^{(n_1)} \mu_{n_2}^{(r')} X^{(n_3)}; ({\bm a},{\bm j}) \right)
&= Y^{(p^r n_1 )} \mu_{\sum_{i=0}^{r-1} p^i b_i +p^r n_2}^{(r+r')} X^{(p^r n_3)} 
B^{({\bm \varepsilon})} ({\bm a},{\bm j}) \\
&= Y^{(p^r n_1 )} \mu_{\sum_{i=0}^{r-1} p^i b_i}^{(r)} 
{\rm Fr'}^r \left( \mu_{n_2'}^{(r')} \right) X^{(p^r n_3)} 
B^{({\bm \varepsilon})} ({\bm a},{\bm j}) \\
&= Y^{(p^r n_1 )}
{\rm Fr'}^r \left( \mu_{n_2'}^{(r')} \right) X^{(p^r n_3)} 
\mu_{\sum_{i=0}^{r-1} p^i b_i}^{(r)} B^{({\bm \varepsilon})} ({\bm a},{\bm j}) \\
&= {\rm Fr'}^r \left( Y^{(n_1)}\mu_{n_2'}^{(r')} X^{(n_3)} \right) 
B^{({\bm \varepsilon})} ({\bm a},{\bm j}),
\end{align*}
and the first claim is proved. It is easy to check the equivalences in the second claim.  
Therefore, the proposition follows. $\square$ \\

\noindent {\bf Remark.} As a special case for $s=t=1$ in (ii), we obtain   
$$Y^{(p^i)}X^{(p^i)}B^{({\bm \varepsilon})}({\bm a},{\bm j}) 
= \gamma_0(a_i,j_i)B^{({\bm \varepsilon})}({\bm a},{\bm j})+4j_i^2
B^{({\bm \varepsilon}+{\bm e}_{i+1})}({\bm a},{\bm j}), $$ 
$$X^{(p^i)}Y^{(p^i)}B^{({\bm \varepsilon})}({\bm a},{\bm j}) 
= \widetilde{\gamma}_0(a_i,j_i)B^{({\bm \varepsilon})}({\bm a},{\bm j})+4j_i^2 
B^{({\bm \varepsilon}+{\bm e}_{i+1})}({\bm a},{\bm j})$$
if $\varepsilon_i =0$ and 
$$Y^{(p^i)}X^{(p^i)}B^{({\bm \varepsilon})}({\bm a},{\bm j}) 
= \gamma_0(a_i,j_i)B^{({\bm \varepsilon})}({\bm a},{\bm j}), $$ 
$$X^{(p^i)}Y^{(p^i)}B^{({\bm \varepsilon})}({\bm a},{\bm j}) 
= \widetilde{\gamma}_0(a_i,j_i)B^{({\bm \varepsilon})}({\bm a},{\bm j})$$
if $\varepsilon_i =1$. 
 \\ 
\

Now we shall give an expression of $B^{({\bm \varepsilon})}({\bm a}, {\bm j})$ for 
$({\bm a},{\bm j}) = \left( (a_i,j_i)\right)_{i=0}^{r-1} \in \mathcal{P}_{\mathbb{Z}}^r$ and 
${\bm \varepsilon}=(\varepsilon_0, \dots, \varepsilon_{r-1}) \in \mathbb{F}_2^r$ which 
generalizes 
Proposition  \ref{PropertiesOfB1} (iii). If we write $B^{(\varepsilon_i)}(a_i,j_i)$ as 
$$B^{(\varepsilon_i)}(a_i,j_i)=\mu_{a_i} \sum_{m_i=n^{(\varepsilon_i)}(a_i,j_i)}^{p-1} 
c^{(\varepsilon_i)}_{m_i}(a_i,j_i) Y^{m_i} X^{m_i} = 
\mu_{a_i} \sum_{\widetilde{m}_i=\widetilde{n}^{(\varepsilon_i)}(a_i,j_i)}^{p-1} 
\widetilde{c}^{(\varepsilon_i)}_{\widetilde{m}_i}(a_i,j_i) X^{\widetilde{m}_i} Y^{\widetilde{m}_i}$$
following the proposition, 
set ${\bm m}= (m_0, \dots, m_{r-1}), 
\widetilde{\bm m}= (\widetilde{m}_0, \dots, \widetilde{m}_{r-1})$, and  
$c^{({\bm \varepsilon})}_{\bm m}({\bm a}, {\bm j})= 
\prod_{i=0}^{r-1}c^{(\varepsilon_i)}_{m_i}(a_i,j_i)$, 
$\widetilde{c}^{({\bm \varepsilon})}_{\widetilde{\bm m}}({\bm a}, {\bm j})= 
\prod_{i=0}^{r-1}\widetilde{c}^{(\varepsilon_i)}_{\widetilde{m}_i}(a_i,j_i)$.  \\
 
\begin{Prop}\label{FormulasOfB2}
For $({\bm a},{\bm j}) = \left( (a_i,j_i)\right)_{i=0}^{r-1} \in \mathcal{P}_{\mathbb{Z}}^r$ and 
${\bm \varepsilon}=(\varepsilon_0, \dots, \varepsilon_{r-1}) \in \mathbb{F}_2^r$, we have 
\begin{align*}
B^{({\bm \varepsilon})}({\bm a}, {\bm j}) 
&= \mu_{\sum_{i=0}^{r-1} p^i b_i}^{(r)} \sum_{\bm m} 
c^{({\bm \varepsilon})}_{\bm m}({\bm a}, {\bm j}) 
\left( \prod_{i=0}^{r-1} Y^{(p^i)m_i} \right) \left( \prod_{i=0}^{r-1} X^{(p^i)m_i} \right) \\
&= \mu_{\sum_{i=0}^{r-1} p^i b_i}^{(r)} \sum_{\widetilde{\bm m}} 
\widetilde{c}^{({\bm \varepsilon})}_{\widetilde{\bm m}}({\bm a}, {\bm j}) 
\left( \prod_{i=0}^{r-1} X^{(p^i)\widetilde{m}_i} \right) 
\left( \prod_{i=0}^{r-1} Y^{(p^i)\widetilde{m}_i} \right),
\end{align*} 
where ${\bm m}= (m_0, \dots, m_{r-1})$ runs through the elements in 
$\mathbb{Z}^r$ satisfying 
$n^{(\varepsilon_i)}(a_i,j_i) \leq m_i \leq p-1$ for $0 \leq i \leq r-1$ and 
$\widetilde{\bm m}= (\widetilde{m}_0, \dots, \widetilde{m}_{r-1})$ runs through the 
elements in $\mathbb{Z}^r$ satisfying 
$\widetilde{n}^{(\varepsilon_i)}(a_i,j_i) \leq \widetilde{m}_i \leq p-1$ for $0 \leq i \leq r-1$. 
\end{Prop}

\noindent {\itshape Proof.} We shall use induction on $r$. It is clear for $r=1$ by Proposition \ref{PropertiesOfB1} (iii), so we may assume that $r \geq 2$. For simplification of notation, set 
$${\bm \varepsilon}' = (\varepsilon_1,\dots, \varepsilon_{r-1}),\ \ 
({\bm a}', {\bm j}')= \left( (a_i,j_i)\right)_{i=1}^{r-1},\ \ 
{\bm m}'= (m_1, \dots, m_{r-1}),$$
$${\bf X}^{\bm m}= \prod_{i=0}^{r-1} X^{(p^i)m_i} , \ \ 
{\bf Y}^{\bm m}= \prod_{i=0}^{r-1} Y^{(p^i)m_i}, \ \ 
{\bf X}^{{\bm m}'}= \prod_{i=1}^{r-1} X^{(p^{i-1})m_i} , \ \ 
{\bf Y}^{{\bm m}'}= \prod_{i=1}^{r-1} Y^{(p^{i-1})m_i},$$ 
$${\bm \mu}=\mu_{\sum_{i=0}^{r-1} p^i b_i}^{(r)},\ \ 
{\bm \mu}'=\mu_{\sum_{i=1}^{r-1} p^{i-1} b_i}^{(r-1)}.$$
Then by induction on $r$ and Propositions \ref{PropertiesOfBZ} (v), 
\ref{PropertiesOfZ} (iii), \ref{PropertiesOfB1} (iii), 
and \ref{CommProp} we have 
\begin{align*}
\lefteqn{B^{({\bm \varepsilon})}({\bm a}, {\bm j}) 
= Z^{(\varepsilon_0)} \left( B^{({\bm \varepsilon}')}({\bm a}', {\bm j}') ;(a_0,j_0) \right) } \\
&= Z^{(\varepsilon_0)} \left( {\bm \mu}' \sum_{{\bm m}'} 
c_{{\bm m}'}^{({\bm \varepsilon}')}({\bm a}',{\bm j}') {\bf Y}^{{\bm m}'} 
{\bf X}^{{\bm m}'};(a_0,j_0) \right) \\
&=  {\bm \mu}  Z^{(\varepsilon_0)} \left(  \sum_{{\bm m}'} 
c_{{\bm m}'}^{({\bm \varepsilon}')}({\bm a}',{\bm j}') {\rm Fr}\left( {\rm Fr}'
\left({\bf Y}^{{\bm m}'}\right) 
{\rm Fr}' \left({\bf X}^{{\bm m}'}\right) \right);(a_0,j_0) \right) \\
&= {\bm \mu}  \sum_{{\bm m}'} 
c_{{\bm m}'}^{({\bm \varepsilon}')}({\bm a}',{\bm j}')  {\rm Fr}'\left({\bf Y}^{{\bm m}'}\right) 
{\rm Fr}'\left({\bf X}^{{\bm m}'}\right) B^{(\varepsilon_0)}(a_0,j_0)  \\
&= {\bm \mu}  \sum_{{\bm m}'} 
c_{{\bm m}'}^{({\bm \varepsilon}')}({\bm a}',{\bm j}')  {\rm Fr}'\left({\bf Y}^{{\bm m}'}\right) 
 B^{(\varepsilon_0)}(a_0,j_0) {\rm Fr}'\left({\bf X}^{{\bm m}'}\right) \\
&= {\bm \mu}  \sum_{{\bm m}'} 
\sum_{m_0=n^{\left(\varepsilon_0\right)}(a_0,j_0)}^{p-1}  c_{m_0}^{(\varepsilon_0)}(a_0,j_0)
c_{{\bm m}'}^{({\bm \varepsilon}')}({\bm a}',{\bm j}')  {\rm Fr}'\left({\bf Y}^{{\bm m}'}\right) 
Y^{m_0}  X^{m_0} {\rm Fr}'\left({\bf X}^{{\bm m}'}\right) \\
&= {\bm \mu}  \sum_{{\bm m}} 
c_{{\bm m}}^{({\bm \varepsilon})}({\bm a},{\bm j})  {\bf Y}^{{\bm m}} 
{\bf X}^{{\bm m}}, 
\end{align*}
and the first equality follows. Similarly, we can show that 
$$B^{({\bm \varepsilon})}({\bm a}, {\bm j}) 
= \mu_{\sum_{i=0}^{r-1} p^i b_i}^{(r)} \sum_{\widetilde{\bm m}} 
\widetilde{c}^{({\bm \varepsilon})}_{\widetilde{\bm m}}({\bm a}, {\bm j}) 
\left( \prod_{i=0}^{r-1} X^{(p^i)\widetilde{m}_i} \right) 
\left( \prod_{i=0}^{r-1} Y^{(p^i)\widetilde{m}_i} \right)$$
and the proposition is proved. $\square$ \\

Now we recall the result on a relation between the elements $B^{({\bm 1})}({\bm a}, {\bm j})$ 
and simple $\mathcal{A}_r$- or $\mathcal{U}_r$-modules. Let 
$({\bm a}, {\bm j})=\left( (a_i,j_i)\right)_{i=0}^{r-1} \in \mathcal{P}_{\mathbb{Z}}^r$. Then 
$\mathcal{A}_r B^{({\bm 1})}({\bm a}, {\bm j})=k B^{({\bm 1})}({\bm a}, {\bm j})$ is a simple 
$\mathcal{A}_r$-module where any nonzero element has $\mathcal{U}_r^0$-weight 
$\sum_{i=0}^{r-1}p^i b_i$ and  is acted on by $Y^{(p^i)}X^{(p^i)}$ as a multiplication 
by $\gamma_0(a_i,j_i)$ for $0 \leq i \leq r-1$, and the set 
$\left\{ \left. \mathcal{A}_r B^{({\bm 1})}({\bm a}, {\bm j})\ \right| \ 
({\bm a}, {\bm j}) \in \mathcal{P}^r \right\}$  forms a complete set of representatives 
of isomorphism classes of simple $\mathcal{A}_r$-modules. On the other hand,   
$\mathcal{U}_r B^{({\bm 1})}({\bm a}, {\bm j})$  is a simple 
$\mathcal{U}_r$-module and 
conversely, any simple $\mathcal{U}_r$-module is isomorphic to 
$\mathcal{U}_r B^{({\bm 1})}({\bm a}, {\bm j})$ for some $({\bm a}, {\bm j}) \in \mathcal{P}_{\mathbb{Z}}^r$ (of course $({\bm a}, {\bm j})$ can be chosen 
in $\mathcal{P}^r$, but it is not determined uniquely 
in general even if it lies there). For details, see 
\cite[\S 3]{yoshii18} and \cite[\S 5]{yoshii22}. \\

\begin{Prop}\label{PropertiesOfB2}
Let $({\bm a}, {\bm j}) = \left( (a_i,j_i)\right)_{i=0}^{r-1} \in \mathcal{P}^r_{\mathbb{Z}}$ and 
${\bm \varepsilon} \in \mathbb{F}^r_2$. Then the following hold. \\ 

\noindent {\rm (i)} $B^{({\bm \varepsilon})}(-{\bm a},{\bm j})$ has 
$\mathcal{U}_r^0$-weight $-\sum_{i=0}^{r-1}p^i b_i$. \\

\noindent {\rm (ii)} $\mathcal{T}_1 \left( B^{({\bm \varepsilon})}({\bm a},{\bm j})\right)=
\mathcal{T}_2 \left( B^{({\bm \varepsilon})}({\bm a},{\bm j})\right)=
B^{({\bm \varepsilon})}(-{\bm a},{\bm j}).$
\end{Prop}
\ 

\noindent {\itshape Proof.} Without loss of generality, we may assume that 
$({\bm a},{\bm j}) \in \mathcal{P}^r$ and ${\bm \varepsilon} \in \mathcal{X}_r({\bm a},{\bm j})$. 

To prove (i), it is enough to check that $-b_i=\iota(-a_i,j_i)$ 
for $0 \leq i \leq r-1$, where $\iota$ is the map defined in Definition \ref{Def2} (1). 

Suppose that $\left(-a_i,j_i \right)$ 
satisfies {\rm (A)} or {\rm (C)}. Then $b_i =a_i$, since 
$(a_i,j_i)$ 
satisfies {\rm (D)} or {\rm (B)}. It follows from the fact 
$a_i \neq 0$ that $1 \leq p-b_i \leq p-1$ and hence that 
$-b_i=(-a_i)\ {\bf mod}\ p-p = \iota(-a_i,j_i)$. 

Suppose that  $\left(-a_i,j_i \right)$ 
satisfies {\rm (B)} or {\rm (D)} and that $a_i \neq 0$.
Then $b_i =a_i- p$, since $(a_i,j_i)$ 
satisfies {\rm (C)} or {\rm (A)}. It follows from the fact 
$a_i \neq 0$ that $1 \leq -b_i \leq p-1$ and hence that 
$-b_i=(-a_i)\ {\bf mod}\ p=\iota(-a_i,j_i)$. 

Finally, suppose that $a_i=0$. Then, since 
$$\left((-a_i)\ {\bf mod}\ p,j_i \right)= (a_i,j_i)= (0,j_i)$$
satisfies {\rm (B)}, we have 
$-b_i = -a_i=(-a_i)\ {\bf mod}\ p=\iota(-a_i,j_i)$, and 
(i) follows.

Now we shall prove (ii). Since $\mathcal{T}_1=\mathcal{T}_2$ on $\mathcal{A}$,  we have 
$\mathcal{T}_1 \left( B^{({\bm \varepsilon})}({\bm a},{\bm j})\right)=
\mathcal{T}_2 \left( B^{({\bm \varepsilon})}({\bm a},{\bm j})\right)$. 
So we only have to show that $\mathcal{T}_1 
\left( B^{({\bm \varepsilon})}({\bm a},{\bm j})\right)=
B^{({\bm \varepsilon})}(-{\bm a},{\bm j})$. Keep in mind that 
${\bm \varepsilon} \in \mathcal{X}_r({\bm a},{\bm j})$. 
 We use induction on the number of nonzero entries of ${\bm \varepsilon}$. 
Suppose that ${\bm \varepsilon} ={\bm 0}$. Since all 
$\mathcal{T}_1 \left(B^{({\bm 0})}({\bm a},{\bm j}) \right)$ with 
$({\bm a},{\bm j}) \in \mathcal{P}^r$ are pairwise orthogonal primitive idempotents 
in the commutative subalgebra $\mathcal{A}_r$ whose 
sum is the unity $1$, we must have 
$\mathcal{T}_1 \left(B^{({\bm 0})}({\bm a},{\bm j}) \right)=
B^{({\bm 0})}\left(\widetilde{\bm a},\widetilde{\bm j}\right) $ for some 
$\left(\widetilde{\bm a},\widetilde{\bm j} \right) = 
\left( \left(\widetilde{a}_i,\widetilde{j}_i\right) \right)_{i=0}^{r-1} 
\in \mathcal{P}_{\mathbb{Z}}^r$ 
(see \cite[ch. 1. Theorem 4.6 (i)]{nagao-tsushimabook}). Thus, to show that 
$\mathcal{T}_1 \left(B^{({\bm 0})}({\bm a},{\bm j}) \right)=
B^{({\bm 0})}\left(-{\bm a},{\bm j}\right)$, it is enough to check that the image of 
$\mathcal{T}_1 \left(B^{({\bm 0})}({\bm a},{\bm j}) \right)$ in the quotient 
$\mathcal{A}_r$-module  
$\mathcal{A}_r / {\rm rad} \mathcal{A}_r$ generates a 
(one-dimensional) simple $\mathcal{A}_r$-module 
isomorphic to $\mathcal{A}_r B^{({\bm 1})}\left(-{\bm a},{\bm j}\right)$.  
By Proposition \ref{PropertiesOfBZ} (and its remark) we have 
\begin{align*}
Y^{(p^i)}X^{(p^i)} \mathcal{T}_1 \left(B^{({\bm 0})}({\bm a},{\bm j}) \right) 
&= \mathcal{T}_1 \left(X^{(p^i)}Y^{(p^i)} B^{({\bm 0})}({\bm a},{\bm j}) \right) \\
&= \widetilde{\gamma}_0(a_i,j_i)  \mathcal{T}_1 \left(B^{({\bm 0})}({\bm a},{\bm j}) \right)
+4j_i^2 \mathcal{T}_1 \left(B^{({\bm e}_{i+1})}({\bm a},{\bm j}) \right) \\
&= {\gamma}_0(-a_i,j_i)  \mathcal{T}_1 \left(B^{({\bm 0})}({\bm a},{\bm j}) \right)
+4j_i^2 \mathcal{T}_1 \left(B^{({\bm e}_{i+1})}({\bm a},{\bm j}) \right)
\end{align*}
for $0 \leq i \leq r-1$ and 
$\mathcal{T}_1 \left(B^{({\bm 0})}({\bm a},{\bm j}) \right)$ 
has the same $\mathcal{U}_r^0$-weight $-\sum_{s=0}^{r-1} p^s b_s$ as 
$B^{({\bm 0})}\left(-{\bm a},{\bm j}\right)$: 
\begin{align*}
{H \choose n} \mathcal{T}_1 \left(B^{({\bm 0})}({\bm a},{\bm j}) \right) 
&= \mathcal{T}_1 \left({-H \choose n} B^{({\bm 0})}({\bm a},{\bm j}) \right) \\
&= {-\sum_{s=0}^{r-1} p^s b_s \choose n} \mathcal{T}_1 \left(B^{({\bm 0})}({\bm a},{\bm j}) \right)
\end{align*}
for $0 \leq n \leq p^r-1$. 
Note that  $4j_i^2\mathcal{T}_1 \left(B^{({\bm e}_{i+1})}({\bm a},{\bm j}) \right)$ lies in 
${\rm rad} \mathcal{A}_r$. Indeed, we have $4j_i^2=0$ in $\mathbb{F}_p$ if 
$(a_i,j_i)$ satisfies (E), and  
$B^{({\bm e}_{i+1})}({\bm a},{\bm j})$ 
(and hence $\mathcal{T}_1 \left(B^{({\bm e}_{i+1})}({\bm a},{\bm j}) \right) $) lies in 
${\rm rad} \mathcal{A}_r$ otherwise (see \cite[Proposition 3.10]{yoshii18}).  
Therefore, the image of the $\mathcal{A}_r$-module 
$\mathcal{A}_r \mathcal{T}_1 \left(B^{({\bm 0})}({\bm a},{\bm j}) \right)$ in   
$\mathcal{A}_r / {\rm rad} \mathcal{A}_r$ is isomorphic to 
$\mathcal{A}_r B^{({\bm 1})}\left(-{\bm a},{\bm j}\right)$ and hence we obtain 
$\mathcal{T}_1 \left(B^{({\bm 0})}({\bm a},{\bm j})\right)=
B^{({\bm 0})}(-{\bm a},{\bm j})$. 

Suppose that ${\bm \varepsilon} =(\varepsilon_0,\dots, \varepsilon_{r-1}) \neq {\bm 0}$. 
There is an integer  $i$ such that  $\varepsilon_i =1$. Then it follows from the assumption 
${\bm \varepsilon} \in \mathcal{X}_r({\bm a},{\bm j})$ that  
$(a_i,j_i)$ does not satisfy (E) and hence 
that $4j_i^2 \neq 0$ in $\mathbb{F}_p$. 
Using induction we obtain 
\begin{align*}
\lefteqn{\mathcal{T}_1 \left(B^{({\bm \varepsilon})}({\bm a},{\bm j}) \right)
=\mathcal{T}_1 \left((4j_i^2)^{-1} \left( Y^{(p^i)}X^{(p^i)} -\gamma_0(a_i,j_i)\right) 
B^{({\bm \varepsilon}-{\bm e}_{i+1})}({\bm a},{\bm j}) \right)} \\
&=  (4j_i^2)^{-1} \left( X^{(p^i)}Y^{(p^i)} -\gamma_0(a_i,j_i)\right) 
\mathcal{T}_1 \left( B^{({\bm \varepsilon}-{\bm e}_{i+1})}({\bm a},{\bm j}) \right) \\
&=  (4j_i^2)^{-1} \left( X^{(p^i)}Y^{(p^i)} -\gamma_0(a_i,j_i)\right) 
B^{({\bm \varepsilon}-{\bm e}_{i+1})}(-{\bm a},{\bm j}) \\
&= (4j_i^2)^{-1} \left( \widetilde{\gamma}_0(-a_i,j_i) -\gamma_0(a_i,j_i)\right) 
B^{({\bm \varepsilon}-{\bm e}_{i+1})}(-{\bm a},{\bm j}) +
 (4j_i^2)^{-1} \cdot 4j_i^2 B^{({\bm \varepsilon})}(-{\bm a},{\bm j}) \\
&= B^{({\bm \varepsilon})}(-{\bm a},{\bm j}),
\end{align*}
as required. $\square$ \\

The following theorem is given in \cite[Corollary 5.5]{yoshii22}, which is used to 
prove the main theorem in this paper. \\

\begin{The}\label{BasisThm}
Let $\mathcal{V} $ denote the subset 
$$
\left\{  
B^{({\bm 0})}\left( ({\bm a}, {\bm j}); 
(t_0, \dots, t_{r-1})\right)\ \left|\ 
\begin{array}{l}
{({\bm a},{\bm j})= \left( (a_i,j_i)\right)_{i=0}^{r-1}\in \mathcal{P}^r,} \\
{-\widetilde{n}^{(0)}(a_i,j_i) \leq t_i \leq n^{(0)}(a_i,j_i),\ \forall i} 
\end{array}
\right. \right\}$$
of  $\bigcup_{({\bm a},{\bm j}) \in \mathcal{P}^r}
\mathcal{B}_r \left(({\bm a}, {\bm j}), {\bm 0}\right)$. Then  
its complement 
$\bigcup_{({\bm a},{\bm j}) \in \mathcal{P}^r}
\mathcal{B}_r \left(({\bm a}, {\bm j}), {\bm 0}\right)
\backslash \mathcal{V}$ 
forms a $k$-basis of  ${\rm rad}\mathcal{U}_r$, and 
the image of $\mathcal{V}$ in  the quotient space  
$\mathcal{U}_r/{\rm rad}\mathcal{U}_r$  forms its $k$-basis. 
\end{The}

\section{Generators of the radical of $\mathcal{U}_r$}
In this section, as a main result  we prove that each of 
certain subsets of  $\mathcal{U}_r$ 
generates ${\rm rad} \mathcal{U}_r$ as an ideal.  
The result improves the main one in \cite{wong83}. If $p$ is odd, set 
$$h(\nu, i)= {H+2p^i \nu \choose 2p^i \nu} +{H+2p^i \nu-1 \choose 2p^i \nu} \in 
\mathcal{U}_{i+1}^0$$
for $\nu \in \{ 1,2,\dots, (p-1)/2\}$ and $i \in \mathbb{Z}_{\geq 0}$. \\

The main result in this paper is as follows. \\

\begin{The}\label{MainThm}
{\rm (i)} Suppose that $p$ is odd. Let $\nu_l$ be  integers with 
$1 \leq \nu_l \leq (p-1)/2$ for $l \in \{0,  \dots, r-1 \}$. Then the  set 
$$\left\{ \left. h(\nu_i,i) X^{(p^i)p-\nu_i}, 
Y^{(p^i)p-\nu_i} h(\nu_i,i)\ \right|\ 0 \leq i \leq r-1 \right\}$$
generates the Jacobson radical ${\rm rad} \mathcal{U}_r$ as  a two-sided ideal 
of $\mathcal{U}_r$.\\ \\
{\rm (ii)} Suppose that  $p=2$. Then the  set 
$$\{ \mu_m^{(i+1)} X^{(m)} X^{(2^i)}, Y^{(2^i)} Y^{(m)} \mu_m^{(i+1)} 
\ |\  0 \leq i \leq r-1, 0 \leq m \leq 2^i-1 \}$$
generates the Jacobson radical ${\rm rad} \mathcal{U}_r$ as  a two-sided  ideal of 
$\mathcal{U}_r$. 
\end{The}

\noindent {\bf Remark.} (a) The generators in (i) for $r=1$ are somewhat different from those 
given by Wong in \cite{wong83}. \\

\noindent (b) The cardinality of the generating set in (i) is $2r$, 
whereas that in (ii) is $2(2^r-1)$. \\

Before  proving the theorem, we shall make some preparations. \\  

\begin{Lem}\label{Lemma1}
Suppose that $p$ is odd. Then we have 
$$\mathcal{T}_1 \left( h(\nu,0) X^{p-\nu} \right)
= (-1)^{\nu+1} Y^{p-\nu} h(\nu,0)$$
in $\mathcal{U}$ for any $\nu \in \{ 1,2, \dots, (p-1)/2\}$.
\end{Lem}

\noindent {\itshape Proof.} It is easy to check that 
$\mathcal{T}_1 \left( {H+2\nu \choose 2\nu} \right) = {H-1 \choose 2\nu}$ and 
$\mathcal{T}_1 \left( {H+2\nu-1 \choose 2\nu} \right) = {H \choose 2\nu}$. 
Then we obtain  
\begin{align*}
\mathcal{T}_1 \left( h(\nu,0) X^{p-\nu}\right)
&= \mathcal{T}_1 \left( \left(  {H+2 \nu \choose 2 \nu} +
{H+2 \nu-1 \choose 2 \nu}\right) X^{p-\nu}\right) \\
&= \left( {H-1 \choose 2\nu}+ {H \choose 2\nu}\right) \cdot 
(-1)^{p-\nu}Y^{p-\nu} \\
&= (-1)^{p-\nu}Y^{p-\nu}\left( {H-2(p-\nu)-1 \choose 2\nu}+ {H-2(p-\nu) \choose 2\nu}\right) \\
&=  (-1)^{\nu+1}Y^{p-\nu}\left( {H+2\nu-1 \choose 2\nu}+ {H+2\nu \choose 2\nu}\right) \\
&= (-1)^{\nu+1}Y^{p-\nu}h(\nu,0), 
\end{align*}
as required. $\square$ 
\ \\ 

\noindent {\bf Remark.} Since $\mathcal{T}_1$ is involutive, we also see that 
$$\mathcal{T}_1 \left( Y^{p-\nu} h(\nu,0)  \right)
= (-1)^{\nu+1} h(\nu,0)X^{p-\nu} $$
for any $\nu \in \{ 1,2, \dots, (p-1)/2\}$ if $p$ is odd. \\

\begin{Lem}\label{Lemma2}
Let $\nu \in \{ 1,2, \dots, (p-1)/2\}$ and $m \in \mathbb{Z}$. Then we have 
$$\mbox{${m+1 \choose 2\nu}+{m \choose 2\nu} \neq 0$ in $\mathbb{F}_p$} 
\Longleftrightarrow  2\nu-1 \leq m\ {\bf mod}\ p \leq p-1.$$ 
\end{Lem}

\noindent {\itshape Proof.} Suppose that $0 \leq m\ {\bf mod}\ p \leq 2\nu-2$. 
Then we have 
${m+1 \choose 2\nu}+ {m \choose 2\nu} = 0+0= 0$ in $\mathbb{F}_p$. 

Conversely, suppose that $2\nu-1 \leq m\ {\bf mod}\ p \leq p-1$. 
Since $0 \leq 2\nu \leq p-1$, without loss of generality we may assume that 
$0 \leq m \leq p-1$. Hence we have $2\nu-1 \leq m \leq p-1$. If 
$m=p-1$, we have 
$${m+1 \choose 2\nu}+ {m \choose 2\nu}=0+{p-1 \choose 2\nu} \neq 0$$
in $\mathbb{F}_p$. If 
$m=2\nu-1$, we have 
$${m+1 \choose 2\nu}+ {m \choose 2\nu}={2\nu \choose 2\nu}+0 =1\neq 0$$
in $\mathbb{F}_p$. Finally, if $2\nu-1< m < p-1$, we have 
$${m+1 \choose 2\nu}+ {m \choose 2\nu}=
\dfrac{m! \cdot 2(m+1-\nu)}{(m+1-2\nu)! (2\nu)!}\neq 0$$
in $\mathbb{F}_p$, as claimed. $\square$ \\

\begin{Lem}\label{Lemma3}
Suppose that $p$ is odd and let $a$ be an integer with $a \not\equiv -1\ ({\rm mod}\ p)$. 
Then $\mu_a X^{p-1}$ lies in the two-sided ideal 
$$\mathcal{U}_1^{\geq 0} h(\nu,0) X^{p-\nu} \mathcal{U}_1^{\geq 0}$$ 
of $\mathcal{U}_1^{\geq 0}$ 
for any $\nu \in \{ 1,2, \dots, (p-1)/2\}$ and $Y^{p-1}\mu_a$ lies in the two-sided ideal  
$$\mathcal{U}_1^{\leq 0} Y^{p-\nu} h(\nu,0)  \mathcal{U}_1^{\leq 0}$$ 
of $\mathcal{U}_1^{\leq 0}$ for any $\nu \in \{ 1,2, \dots, (p-1)/2\}$. 
\end{Lem}

\noindent {\itshape Proof.} Note that all the elements 
$$\mu_a X^{\lambda} h(\nu,0) X^{p-(\lambda+1)} = 
\left( {a-2\lambda+2\nu \choose 2\nu}+ 
{a-2\lambda+2\nu-1 \choose 2\nu}\right)\mu_a X^{p-1}$$
with $0 \leq \lambda \leq \nu-1$ lie in  
$\mathcal{U}_1^{\geq 0} h(\nu,0) X^{p-\nu} \mathcal{U}_1^{\geq 0}$. 
Since $a \not\equiv -1\ ({\rm mod}\ p)$, the integer $(a+2\nu-1)\ {\bf mod}\ p$ lies in 
 $\{0,\dots,p-1\} \backslash \{2\nu-2\}$. 
If $$2\nu-1 \leq (a+2\nu-1)\ {\bf mod}\ p \leq p-1,$$ 
taking  $\lambda=0$ we have 
${a-2\lambda+2\nu \choose 2\nu} + {a-2\lambda+2\nu-1 \choose 2\nu}\neq 0$ 
in $\mathbb{F}_p$  by Lemma \ref{Lemma2}. 
Suppose that $$0 \leq (a+2\nu-1)\ {\bf mod}\ p \leq 2\nu-3.$$ If the integer 
$(a+2\nu-1)\ {\bf mod}\ p$ is even, taking 
$$\lambda = \dfrac{(a+2\nu-1)\ {\bf mod}\ p}{2}+1\ (\leq \nu-1)$$ 
we have  
${a-2\lambda+2\nu \choose 2\nu} + {a-2\lambda+2\nu-1 \choose 2\nu}
= {p-1 \choose 2\nu}+ {p-2 \choose 2\nu}\neq 0$ 
in $\mathbb{F}_p$ by Lemma \ref{Lemma2}.   
If the integer 
$(a+2\nu-1)\ {\bf mod}\ p$ is odd, taking 
$$\lambda = \dfrac{ \left( (a+2\nu-1)\ {\bf mod}\ p\right)+1}{2}\ (\leq \nu-1)$$ 
we have 
${a-2\lambda+2\nu \choose 2\nu} + {a-2\lambda+2\nu-1 \choose 2\nu}
= 0+{p-1 \choose 2\nu}\neq 0$ in $\mathbb{F}_p$.  
Therefore, we have shown that  
${a-2\lambda+2\nu \choose 2\nu}+{a-2\lambda+2\nu-1 \choose 2\nu} \neq 0$ 
in $\mathbb{F}_p$ for some 
integer $\lambda$ with $0 \leq \lambda \leq \nu-1$ and then 
the element 
$$
\mu_a X^{p-1} 
= \left( {a-2\lambda+2\nu \choose 2\nu}+ 
{a-2\lambda+2\nu-1 \choose 2\nu}\right)^{-1} 
\mu_a X^{\lambda} h(\nu,0) X^{p-(\lambda+1)} 
$$
lies in $\mathcal{U}_1^{\geq 0} h(\nu,0) X^{p-\nu} \mathcal{U}_1^{\geq 0}$. 
Moreover, applying $\mathcal{T}_1$ to the element $\mu_{-a-2} X^{p-1}$ which 
lies in $\mathcal{U}_1^{\geq 0} h(\nu,0) X^{p-\nu} \mathcal{U}_1^{\geq 0}$ 
(note that $-a-2 \not\equiv -1\ ({\rm mod}\ p)$), 
we also see that  
$Y^{p-1} \mu_a$ lies in 
$\mathcal{U}_1^{\leq 0} Y^{p-\nu} h(\nu,0)  \mathcal{U}_1^{\leq 0}$ 
by Lemma \ref{Lemma1}, since 
$$\mathcal{T}_1(\mu_{-a-2} X^{p-1})= \mu_{a+2} \cdot (-1)^{p-1}Y^{p-1}=Y^{p-1}\mu_a$$
by Proposition \ref{FormulasOfMu} (i).  Thus,  the lemma follows. $\square$ \\

\begin{Lem}\label{Lemma4}
Let $(a,j) \in \mathcal{P}_{\mathbb{Z}}$.  Then 
$B^{(1)}(a,j)$ lies in the two-sided ideal 
$$\mathcal{U}_1 \mu_{a+2n^{(0)}(a,j)} X^{p-1} \mathcal{U}_1$$ 
of $\mathcal{U}_1$. 
In particular, if $p$ is odd and $j \neq 0$, then  
$B^{(1)}(a,j)$ lies in the two-sided ideal 
$$\mathcal{U}_1 
h(\nu,0) X^{p-\nu} \mathcal{U}_1$$
of $\mathcal{U}_1$ 
for any 
$\nu \in \{1,2,\dots, (p-1)/2 \}$. 
\end{Lem}

\noindent {\itshape Proof.} We know that the 
element $$Y^{n^{(0)}(a,j)}X^{n^{(0)}(a,j)}B^{(1)}(a,j)= \beta_{n^{(0)}(a,j)}(a,j) B^{(1)}(a,j)$$
is a nonzero scalar multiple of $B^{(1)}(a,j)$ (see Remark (a) of Definition 4.1 in 
\cite{yoshii22}).  Note that $B^{(1)}(a,j) = 
\mu_a X^{\widetilde{n}^{(1)}(a,j)} y$ for some 
nonzero $y \in \mathcal{U}_1$ by 
Proposition \ref{PropertiesOfB1} (iii) and that 
$$X^{n^{(0)}(a,j)}B^{(1)}(a,j)= X^{n^{(0)}(a,j)} \cdot \mu_a 
X^{\widetilde{n}^{(1)}(a,j)} y = \mu_{a+2n^{(0)}(a,j)} X^{p-1}y.$$
Then we have 
$$B^{(1)}(a,j) = \beta_{n^{(0)}(a,j)}(a,j)^{-1} 
Y^{n^{(0)}(a,j)}\mu_{a+2n^{(0)}(a,j)} X^{p-1} y$$
and hence  the first claim is proved. Moreover, if $p$ is odd and $j \neq 0$, 
then Lemma \ref{Lemma3} shows that 
$B^{(1)}(a,j)$ must lie in $\mathcal{U}_1 h(\nu,0) X^{p-\nu} \mathcal{U}_1$ for any 
$\nu \in \{1,2,\dots, (p-1)/2 \}$, since 
$a+2n^{(0)}(a,j) \not\equiv -1\ ({\rm mod}\ p)$. 
Therefore, the lemma follows. $\square$ \\

For $(a,j) \in \mathcal{P}_{\mathbb{Z}}$, recall that $n^{(0)}(a,j) \leq n^{(1)}(a,j)$ and 
$\widetilde{n}^{(0)}(a,j) \leq \widetilde{n}^{(1)}(a,j)$ and that the equalities hold if and only if 
$(a,j)$ satisfies {\rm (E)}. \\

\begin{Lem}\label{Lemma5}
Let $(a,j) \in \mathcal{P}_{\mathbb{Z}}$.  Let $t$ be an integer satisfying 
$n^{(0)}(a,j) < t \leq n^{(1)}(a,j)$ or $-\widetilde{n}^{(1)}(a,j) \leq t < 
-\widetilde{n}^{(0)}(a,j)$ (these occur only if $(a,j)$ does not satisfy {\rm (E)}). 
Then  $B^{(0)}\left( (a,j);t\right)$ lies in the two-sided ideal 
$$
\left\{ \begin{array}{ll}
{\mathcal{U}_1 \mu_{a+2n^{(1)}(a,j)}X^{p-1} \mathcal{U}_1} & 
\mbox{if $n^{(0)}(a,j) < t \leq n^{(1)}(a,j)$,} \\
{\mathcal{U}_1 Y^{p-1}\mu_{a+2n^{(0)}(a,j)} \mathcal{U}_1} &
\mbox{if $-\widetilde{n}^{(1)}(a,j) \leq t < -\widetilde{n}^{(0)}(a,j)$}
\end{array} \right.
$$
of $\mathcal{U}_1$. 
In particular, if $p$ is odd, then  $B^{(0)}\left( (a,j);t\right)$ lies in the two-sided ideal  
$$
\left\{ \begin{array}{ll}
{\mathcal{U}_1 h(\nu,0) X^{p-\nu}\mathcal{U}_1} & 
\mbox{if $n^{(0)}(a,j) < t \leq n^{(1)}(a,j)$,} \\
{\mathcal{U}_1 Y^{p-\nu} h(\nu,0) \mathcal{U}_1} &
\mbox{if $-\widetilde{n}^{(1)}(a,j) \leq t < -\widetilde{n}^{(0)}(a,j)$}
\end{array} \right.
$$
of $\mathcal{U}_1$ for any $\nu \in \{1,2,\dots, (p-1)/2 \}$. 
\end{Lem}

\noindent {\itshape Proof.} Suppose  that 
$n^{(0)}(a,j) < t \leq n^{(1)}(a,j)$. 
Since 
$B^{(1)}\left( (a,j);t\right)=0$ by Proposition \ref{PropertiesOfB1} (iii) and 
Remark (c) of Definition \ref{Def1}, 
it follows from \cite[Lemma 4.6 (i)]{yoshii22} 
that $Y^{n^{(1)}(a,j)-t}X^{n^{(1)}(a,j)}B^{(0)}(a,j)= c B^{(0)}\left( (a,j);t\right)$ 
for some $c \in k$. But the scalar $c$ must be nonzero since 
$$Y^{n^{(1)}(a,j)}X^{n^{(1)}(a,j)}B^{(0)}(a,j)= 
4j^2 \left( \prod_{i=0,\ i \neq {n}^{(0)}(a,j)}^{{n}^{(1)}(a,j)-1}\gamma_i(a,j) \right) B^{(1)}(a,j)
\neq 0.$$
Note that  $B^{(0)}(a,j)= \mu_a X^{\widetilde{n}^{(0)}(a,j)} y$ for some 
nonzero $y \in \mathcal{U}_1$ by 
Proposition \ref{PropertiesOfB1} (iii) and that 
$$X^{n^{(1)}(a,j)}B^{(0)}(a,j)=X^{n^{(1)}(a,j)} \cdot \mu_a X^{\widetilde{n}^{(0)}(a,j)} y= 
\mu_{a+2n^{(1)}(a,j)} X^{p-1}y.$$
Hence the element 
$$B^{(0)}\left( (a,j);t \right)=c^{-1} Y^{n^{(1)}(a,j)-t}
\mu_{a+2n^{(1)}(a,j)} X^{p-1}y$$
 must lie in $\mathcal{U}_1 \mu_{a+2n^{(1)}(a,j)}X^{p-1} \mathcal{U}_1$. 

Suppose that $-\widetilde{n}^{(1)}(a,j) \leq t < -\widetilde{n}^{(0)}(a,j)$. Then the result in the 
last paragraph shows that $B^{(0)}\left( (-a,j);-t \right)$  lies in 
$\mathcal{U}_1 \mu_{-a+2n^{(1)}(-a,j)}X^{p-1} \mathcal{U}_1$, since 
$n^{(0)}(-a,j) < -t \leq n^{(1)}(-a,j)$. Then the element 
$\mathcal{T}_1 \left( B^{(0)}\left( (-a,j);-t \right)\right)= 
(-1)^{-t} B^{(0)}\left( (a,j);t \right)$ lies in 
$\mathcal{T}_1 \left(\mathcal{U}_1 \mu_{-a+2n^{(1)}(-a,j)}X^{p-1} \mathcal{U}_1\right)$. 
Since 
\begin{align*}
\mathcal{T}_1 \left(\mu_{-a+2n^{(1)}(-a,j)}X^{p-1}\right) 
&= (-1)^{p-1} \mu_{a-2n^{(1)}(-a,j)}Y^{p-1} \\
&= \mu_{a-2\widetilde{n}^{(1)}(a,j)}Y^{p-1} \\
&= Y^{p-1} \mu_{a+2(p-1-\widetilde{n}^{(1)}(a,j))} \\
&= Y^{p-1} \mu_{a+2{n}^{(0)}(a,j)} 
\end{align*} 
by Proposition \ref{FormulasOfMu} (i) and Remark (c) of Definition \ref{Def1},  
$ B^{(0)}\left( (a,j);t \right)$ must lie in 
$\mathcal{U}_1 Y^{p-1}\mu_{a+2n^{(0)}(a,j)} \mathcal{U}_1$. 
Therefore, the first claim is proved. 

Finally, suppose that $p$ is odd. Since 
$(a,j)$ does not satisfy (E) (i.e. $j \neq 0$), neither $a+2n^{(0)}(a,j)$ nor 
$a+2n^{(1)}(a,j)$ is congruent to $-1$ modulo $p$. Now the second claim in the lemma 
 follows from the first 
claim and Lemma \ref{Lemma3}. Therefore,  the lemma follows. $\square$ \\

\begin{Lem}\label{Lemma6}
{\rm (i)} Suppose that $p$ is odd. Let $\nu$ be an integer with 
$1 \leq \nu \leq (p-1)/2$. Then   $h(\nu,0)X^{p-\nu} $ and 
$Y^{p-\nu}  h(\nu,0) $ lie in ${\rm rad} \mathcal{U}_1$. \\

\noindent {\rm (ii)} Suppose that $p=2$. Then 
  $\mu_0 X$ and 
$Y  \mu_0$ lie in ${\rm rad} \mathcal{U}_1$. 
\end{Lem}

\noindent {\itshape Proof.} Suppose that $p=2$. Since 
\begin{align*}
{\rm soc}_{\mathcal{U}_1} \mathcal{U}_1 
&= \mathcal{U}_1 B^{(1)}\left( 0,\dfrac{1}{2}\right)
+\mathcal{U}_1 B^{(1)}(1,0)+\mathcal{U}_1 B^{(1)}(1,1) \\
&= \mathcal{U}_1 \mu_0XY + \mathcal{U}_1 \mu_1YX+   \mathcal{U}_1 \mu_1XY \\
&=  k \mu_0XY + k \mu_1YX+ k\mu_1X +k \mu_1XY+k\mu_1Y,
\end{align*}
we easily see that the elements $\mu_0X$ and $Y\mu_0$ annihilate 
${\rm soc}_{\mathcal{U}_1} \mathcal{U}_1$ and hence lie in 
${\rm rad} \mathcal{U}_1$. Thus (ii) follows. 

From now on, suppose that $p$ is odd. Fix $\nu \in \{ 1,2,\dots, (p-1)/2\}$. 
Recall from the remark of Theorem 5.3 in 
\cite{yoshii22} that any simple $\mathcal{U}_1$-module can be written as 
$\mathcal{U}_1 B^{(1)}(a,j)$ for some $(a,j) \in \mathcal{P}_{\mathbb{Z}}$ 
and from Proposition 5.2 (ii) in \cite{yoshii22} that the set 
$$\left\{ \left. B^{(1)}\left( (a,j); t\right)\ \right|\ 
-\widetilde{n}^{(0)}(a,j) \leq t \leq n^{(0)}(a,j)
\right\}$$
forms a $k$-basis of $\mathcal{U}_1 B^{(1)}(a,j)$. 
Fix a pair $(a,j) \in \mathcal{P}_{\mathbb{Z}}$ and an integer 
$t$ with $-\widetilde{n}^{(0)}(a,j) \leq t \leq n^{(0)}(a,j)$. We know that 
$X^{p-\nu} B^{(1)}\left( (a,j); t\right)$ is a nonzero scalar multiple of 
$$B^{(1)}\left( (a,j); t+p-\nu\right),$$ which is not zero only if $t+p-\nu \leq n^{(0)}(a,j)$ 
(see Lemma 4.7 and Proposition 4.5 in \cite{yoshii22}). 
Of course $B^{(1)}(a,j)$ has $\mathcal{U}_1^0$-weight $b$ defined in 
Definition \ref{Def2} (1). 
Suppose that $X^{p-\nu} B^{(1)}\left( (a,j);t\right) \neq 0$ 
(so $t$ satisfies $-\widetilde{n}^{(0)}(a,j) \leq t \leq n^{(0)}(a,j)+\nu-p$). Note that 
\begin{align*}
\lefteqn{h(\nu,0) B^{(1)} \left( (a,j);t+p-\nu \right)} \\
&= \left( {b+2t+p \choose 2\nu}+ {b+2t+p-1 \choose 2\nu}
 \right) B^{(1)} \left( (a,j);t+p-\nu \right).
\end{align*}
Since 
$$b-2 \widetilde{n}^{(0)}(a,j)-1 \leq b+2t-1 \leq b+
2\left( n^{(0)}(a,j) +\nu-p\right)-1,$$
we have $0 \leq b+2t+p-1 \leq 2\nu-2$ by the remark of Definition \ref{Def2} and hence 
${b+2t+p \choose 2\nu}={b+2t+p-1 \choose 2\nu}=0$ in 
$\mathbb{F}_p$. Therefore, we obtain 
$$h(\nu,0) X^{p-\nu} B^{(1)}\left( (a,j); t\right) =0.$$
We also note from this result that 
$$h(\nu,0) X^{p-\nu} B^{(1)}\left( (-a,j); -t\right) =0,$$
since $-\widetilde{n}^{(0)}(-a,j)\leq -t \leq n^{(0)}(-a,j)$. Applying the map $\mathcal{T}_1$ to 
the equality and using Lemma \ref{Lemma1} and Proposition \ref{PropertiesOfB2} 
(ii) we obtain
$$Y^{p-\nu} h(\nu,0)B^{(1)}\left( (a,j); t\right) =0.$$
Therefore, the elements  $h(\nu,0) X^{p-\nu} $ and  $Y^{p-\nu} h(\nu,0)$ 
 annihilate each simple $\mathcal{U}_1$-module and 
hence lie in ${\rm rad} \mathcal{U}_1$, and 
 (i) is proved. $\square$ \\

\begin{Lem}\label{Lemma7}  
The following hold. \\ \\
{\rm (i)} Suppose that $p$ is odd. Let 
$({\bm a}, {\bm j})=\left( (a_i,j_i) \right)_{i=0}^{r-1} \in \mathcal{P}_{\mathbb{Z}}^r$ and  
let $a'$ be an integer 
with $a' \not\equiv -1\ ({\rm mod}\ p)$. Then the element 
$Z^{({\bm 0})}\left(\mu_{a'}X^{p-1}; ({\bm a},{\bm j})\right)$ lies in the two-sided ideal  
$$\mathcal{U}_{r+1} h(\nu,r)
X^{(p^r)p-\nu} \mathcal{U}_{r+1}$$ 
of $\mathcal{U}_{r+1}$ for any $\nu \in \{ 1,2,\dots, (p-1)/2\}$ and 
the element 
$Z^{({\bm 0})}\left(Y^{p-1}\mu_{a'}; ({\bm a},{\bm j})\right)$ lies in the two-sided ideal 
$$\mathcal{U}_{r+1} Y^{(p^r)p-\nu} h(\nu,r) \mathcal{U}_{r+1}$$
 of $\mathcal{U}_{r+1}$  
for any $\nu \in \{ 1,2,\dots, (p-1)/2\}$. \ \\ 

\noindent {\rm (ii)} Suppose that $p=2$. Let $(a_0, j_0) \in \mathcal{P}_{\mathbb{Z}}$ and  
let $a'$ be an integer with 
$0 \leq a' \leq 2^{r-1}-1$. Then the element 
$Z^{(0)}\left(\mu^{(r)}_{a'}X^{(a')}X^{(2^{r-1})}; (a_0,j_0)\right)$ lies in the two-sided ideal 
$$\mathcal{U}_{r+1} \mu^{(r+1)}_{(a_0\ {\bf mod}\ 2)+2a'} 
X^{((a_0\ {\bf mod}\ 2)+2a')} X^{(2^r)}\mathcal{U}_{r+1}$$ of 
$\mathcal{U}_{r+1}$ and 
the element 
$Z^{(0)}\left(Y^{(2^{r-1})}Y^{(a')}\mu^{(r)}_{a'}; (a_0,j_0)\right)$ lies in the two-sided ideal 
$$\mathcal{U}_{r+1} Y^{(2^r)}
Y^{((a_0\ {\bf mod}\ 2)+2a')}\mu^{(r+1)}_{(a_0\ {\bf mod}\ 2)+2a'}  \mathcal{U}_{r+1}$$ 
of $\mathcal{U}_{r+1}$.
\end{Lem}

\noindent {\itshape Proof.} 
Suppose first that $p=2$. Recall that 
$$B^{(0)}\left(2i, \dfrac{1}{2}\right)=\mu_0,\ \ \ 
B^{(0)}(1+2i,0)=\mu_1 YX,\ \ \ B^{(0)}(1+2i,1)=\mu_1 XY$$
for $i \in \mathbb{Z}$. 
Using Propositions \ref{PropertiesOfBZ} (iv), (v), \ref{FormulasOfMu} (i), (ii), 
and \ref{CommProp} we easily see that 
$$
Z^{(0)}\left( \mu_{a'}^{(r)} X^{(a')} X^{(2^{r-1})} ; (a_0,j_0)\right) 
= \mu_{2a'}^{(r+1)} X^{(2a')} X^{(2^r)}, 
$$
$$Z^{(0)}\left(Y^{(2^{r-1})}Y^{(a')}\mu_{a'}^{(r)}; (a_0,j_0)\right) = 
Y^{(2^{r})}Y^{(2a')}\mu_{2a'}^{(r+1)}$$ 
if $(a_0\ {\bf mod}\ 2,j_0)=(0,1/2)$,
$$
Z^{(0)}\left( \mu_{a'}^{(r)} X^{(a')} X^{(2^{r-1})} ; (a_0,j_0)\right) 
= Y\mu_{1+2a'}^{(r+1)} X^{(1+2a')} X^{(2^r)}, 
$$
$$Z^{(0)}\left(Y^{(2^{r-1})}Y^{(a')}\mu_{a'}^{(r)}; (a_0,j_0)\right) = 
Y^{(2^{r})}Y^{(1+2a')}\mu_{1+2a'}^{(r+1)}X$$ 
if $(a_0\ {\bf mod}\ 2,j_0)=(1,0)$, and  
$$
Z^{(0)}\left( \mu_{a'}^{(r)} X^{(a')} X^{(2^{r-1})} ; (a_0,j_0)\right) 
= \mu_{1+2a'}^{(r+1)} X^{(1+2a')} X^{(2^r)}Y, 
$$
$$Z^{(0)}\left(Y^{(2^{r-1})}Y^{(a')}\mu_{a'}^{(r)}; (a_0,j_0)\right) = 
XY^{(2^{r})}Y^{(1+2a')}\mu_{1+2a'}^{(r+1)}$$ 
if $(a_0\ {\bf mod}\ 2,j_0)=(1,1)$. Thus (ii) follows. 

From now on, assume that $p$ is odd. Fix $\nu \in \{1,2,\dots, (p-1)/2\}$. 
We proceed in steps. \\

\noindent {\bf Step 1.} For any integer $\lambda$ with $0 \leq \lambda \leq \nu-1$ 
and any $d \in \mathbb{Z}$, 
each of the two elements 
$$\mu_{\pm \sum_{i=0}^{r-1}p^i\left(b_i+2n^{(0)}(a_i,j_i) \right)+p^r d -2p^r \lambda}^{(r+1)}
 X^{(p^r)p-\nu}$$
is a factor  in a factorization of $Z^{({\bm 0})} \left( \mu_{d} X^{p-1} ; ({\bm a}, {\bm j})\right)$ in $\mathcal{U}_{r+1}$ and 
each of the two elements 
$$Y^{(p^r)p-\nu}\mu_{\pm \sum_{i=0}^{r-1}p^i \left(b_i+2n^{(0)}(a_i,j_i) \right)+
p^r d -2p^r \lambda}^{(r+1)} $$
is a factor  in a factorization of  
$Z^{({\bm 0})} \left( Y^{p-1}\mu_{d}  ; ({\bm a}, {\bm j})\right)$ 
in $\mathcal{U}_{r+1}$. \\

By Proposition \ref{FormulasOfB2} $B^{({\bm 0})}({\bm a}, {\bm j})$ can be written as 
$$
B^{({\bm 0})}({\bm a}, {\bm j}) 
= \mu_{\sum_{i=0}^{r-1} p^i b_i}^{(r)} \sum_{\bm m} 
c^{({\bm 0})}_{\bm m}({\bm a}, {\bm j}) 
\left( \prod_{i=0}^{r-1} Y^{(p^i)m_i} \right) \left( \prod_{i=0}^{r-1} X^{(p^i)m_i} \right),  
$$ 
where ${\bm m}= (m_0, \dots, m_{r-1})$ runs through the elements in $\mathbb{Z}^r$ satisfying 
$n^{(0)}(a_i,j_i) \leq m_i \leq p-1$ for all integers $i$ with $0 \leq i \leq r-1$.  
For simplification of notation, set 
$$n=\sum_{i=0}^{r-1}p^i b_i,\ \ n'=\sum_{i=0}^{r-1} p^i\left(b_i+2n^{(0)}(a_i,j_i) \right),$$
$${\bf X}=\prod_{i=0}^{r-1} X^{(p^i)n^{(0)}(a_i,j_i)},\ \ 
{\bf X}^{\bm m}=\prod_{i=0}^{r-1} X^{(p^i)m_i}, \ \ 
{\bf X}_0^{\bm m}=\prod_{i=0}^{r-1} X^{(p^i)m_i-n^{(0)}(a_i,j_i)},$$
$${\bf Y}=\prod_{i=0}^{r-1} Y^{(p^i)n^{(0)}(a_i,j_i)},\ \ 
{\bf Y}^{\bm m}=\prod_{i=0}^{r-1} Y^{(p^i)m_i}, \ \ 
{\bf Y}_0^{\bm m}=\prod_{i=0}^{r-1} Y^{(p^i)m_i-n^{(0)}(a_i,j_i)}.$$
Note that ${\bf X}^{\bm m}={\bf X}_0^{\bm m}{\bf X}$ and 
${\bf Y}^{\bm m}={\bf Y}{\bf Y}_0^{\bm m}$. 
By Propositions \ref{PropertiesOfBZ} (iv), (v), \ref{FormulasOfMu} (i), (ii), 
and \ref{CommProp}  we have 
\begin{align*}
Z^{({\bm 0})} \left( \mu_{d} X^{p-1} ; ({\bm a}, {\bm j})\right)
&= \mu_{\sum_{i=0}^{r-1}p^i b_i +p^r d}^{(r+1)} X^{(p^r)p-1} B^{({\bm 0})}({\bm a}, {\bm j}) \\
&= \mu_{n+p^r d}^{(r+1)} \mu_n^{(r)}\sum_{\bm m} c^{({\bm 0})}_{\bm m}({\bm a}, {\bm j})  
{\bf Y}^{\bm m} {\bf X}^{\bm m} X^{(p^r)p-1} \\
&= \mu_{n+p^r d}^{(r+1)} \sum_{\bm m} c^{({\bm 0})}_{\bm m}({\bm a}, {\bm j})  
{\bf Y}^{\bm m} {\bf X}^{\bm m} X^{(p^r)p-1} \\
&= {\bf Y} \mu_{n'+p^r d}^{(r+1)} \sum_{\bm m} c^{({\bm 0})}_{\bm m}({\bm a}, {\bm j})  
{\bf Y}_0^{\bm m} {\bf X}_0^{\bm m} X^{(p^r)p-1} {\bf X} \\
&= {\bf Y} \mu_{n'+p^r d}^{(r+1)} X^{(p^r)p-1}
\sum_{\bm m} c^{({\bm 0})}_{\bm m}({\bm a}, {\bm j})  
{\bf Y}_0^{\bm m} {\bf X}_0^{\bm m}  {\bf X} \\
&= {\bf Y} X^{(p^r)\lambda} \mu_{n'+p^r d-2p^r\lambda}^{(r+1)} X^{(p^r)p-(\lambda+1)}
\sum_{\bm m} c^{({\bm 0})}_{\bm m}({\bm a}, {\bm j})  
{\bf Y}_0^{\bm m} {\bf X}^{\bm m}
\end{align*}
for any integer $\lambda$ with $0 \leq \lambda \leq \nu-1$, and hence  
we see that $\mu_{n'+p^r d-2p^r\lambda}^{(r+1)} X^{(p^r)p-\nu}$ is a factor in a 
factorization of 
$Z^{({\bm 0})} \left( \mu_{d} X^{p-1} ; ({\bm a}, {\bm j})\right)$ in $\mathcal{U}_{r+1}$. 

On the other hand, by Proposition \ref{FormulasOfB2} 
again $B^{({\bm 0})}({\bm a}, {\bm j})$ can be also 
written as 
$$
B^{({\bm 0})}({\bm a}, {\bm j}) 
= \mu_{\sum_{i=0}^{r-1} p^i b_i}^{(r)} \sum_{\widetilde{\bm m}} 
\widetilde{c}^{({\bm 0})}_{\widetilde{\bm m}}({\bm a}, {\bm j}) 
\left( \prod_{i=0}^{r-1} X^{(p^i)\widetilde{m}_i} \right) 
\left( \prod_{i=0}^{r-1} Y^{(p^i)\widetilde{m}_i} \right),  
$$ 
where $\widetilde{\bm m}= (\widetilde{m}_0, \dots, \widetilde{m}_{r-1})$ 
runs through the elements in $\mathbb{Z}^r$ satisfying 
$\widetilde{n}^{(0)}(a_i,j_i) \leq \widetilde{m}_i \leq p-1$ for all 
integers $i$ with $0 \leq i \leq r-1$.  For simplification of notation, set 
$$\widetilde{n}'=\sum_{i=0}^{r-1} p^i \left(b_i-2\widetilde{n}^{(0)}(a_i,j_i)\right),$$
$$\widetilde{\bf X}=\prod_{i=0}^{r-1} X^{(p^i)\widetilde{n}^{(0)}(a_i,j_i)},\ \ 
\widetilde{\bf X}^{\widetilde{\bm m}}=\prod_{i=0}^{r-1} X^{(p^i)\widetilde{m}_i}, \ \ 
\widetilde{\bf X}_0^{\widetilde{\bm m}}=\prod_{i=0}^{r-1} X^{(p^i)\widetilde{m}_i
-\widetilde{n}^{(0)}(a_i,j_i)},$$
$$
\widetilde{\bf Y}^{\widetilde{\bm m}}=\prod_{i=0}^{r-1} Y^{(p^i)\widetilde{m}_i}.$$
Note that $\widetilde{\bf X}^{\widetilde{\bm m}}=\widetilde{\bf X}
\widetilde{\bf X}_0^{\widetilde{\bm m}} $ and that 
$\widetilde{n}'=-n'$, since 
$b_i -2\widetilde{n}^{(0)}(a_i,j_i) = - \left( b_i +2n^{(0)}(a_i,j_i) \right)$ 
for any $i$ with $0 \leq i \leq r-1$ (see the remark of Definition \ref{Def2}). 
By Propositions \ref{PropertiesOfBZ} (iv), (v), \ref{FormulasOfMu} (i), (ii), 
and \ref{CommProp} we have 
\begin{align*}
Z^{({\bm 0})} \left( \mu_{d} X^{p-1} ; ({\bm a}, {\bm j})\right)
&= \mu_{\sum_{i=0}^{r-1}p^i b_i +p^r d}^{(r+1)} X^{(p^r)p-1} B^{({\bm 0})}({\bm a}, {\bm j}) \\
&= \mu_{n+p^r d}^{(r+1)} X^{(p^r)p-1} \mu_n^{(r)} \sum_{\widetilde{\bm m}} 
\widetilde{c}^{({\bm 0})}_{\widetilde{\bm m}}({\bm a}, {\bm j})  
\widetilde{\bf X}^{\widetilde{\bm m}} \widetilde{\bf Y}^{\widetilde{\bm m}}  \\
&= \mu_{n+p^r d}^{(r+1)} \mu_n^{(r)} X^{(p^r)p-1} \sum_{\widetilde{\bm m}} 
\widetilde{c}^{({\bm 0})}_{\widetilde{\bm m}}({\bm a}, {\bm j})  
\widetilde{\bf X}^{\widetilde{\bm m}} \widetilde{\bf Y}^{\widetilde{\bm m}}  \\
&= \widetilde{\bf X} \mu_{\widetilde{n}'+p^r d}^{(r+1)} X^{(p^r)p-1}
\sum_{\widetilde{\bm m}} \widetilde{c}^{({\bm 0})}_{\widetilde{\bm m}}({\bm a}, {\bm j})  
\widetilde{\bf X}_0^{\widetilde{\bm m}} 
\widetilde{\bf Y}^{\widetilde{\bm m}}   \\
&= \widetilde{\bf X} \mu_{-{n}'+p^r d}^{(r+1)} X^{(p^r)p-1}
\sum_{\widetilde{\bm m}} \widetilde{c}^{({\bm 0})}_{\widetilde{\bm m}}({\bm a}, {\bm j})  
\widetilde{\bf X}_0^{\widetilde{\bm m}} 
\widetilde{\bf Y}^{\widetilde{\bm m}} \\
&= \widetilde{\bf X} X^{(p^r)\lambda} \mu_{-n'+p^r d-2p^r \lambda}^{(r+1)} X^{(p^r)p-(\lambda+1)}
\sum_{\widetilde{\bm m}} \widetilde{c}^{({\bm 0})}_{\widetilde{\bm m}}({\bm a}, {\bm j})  
\widetilde{\bf X}_0^{\widetilde{\bm m}} 
\widetilde{\bf Y}^{\widetilde{\bm m}}   
\end{align*}
for any integer $\lambda$ with $0 \leq \lambda \leq \nu-1$, and hence  
we see that $\mu_{-n'+p^r d-2p^r\lambda}^{(r+1)} X^{(p^r)p-\nu}$ is a  factor in a 
factorization of 
$Z^{({\bm 0})} \left( \mu_{d} X^{p-1} ; ({\bm a}, {\bm j})\right)$ 
in $\mathcal{U}_{r+1}$. 
Therefore, the first claim in the step follows. 

We shall show the second claim. Since 
$$Y^{(p^r)p-1}\mu_{\pm n'+p^r d}^{(r+1)}  = 
Y^{(p^r)p-(\lambda+1)}\mu_{\pm n'+p^r d-2p^r\lambda}^{(r+1)} Y^{(p^r)\lambda} $$
for any integer $\lambda$ with $0 \leq \lambda \leq \nu-1$, it is enough to check that 
each of the elements $Y^{(p^r)p-1}\mu_{\pm n'+p^r d}^{(r+1)}$ is a factor in 
a factorization of 
$Z^{({\bm 0})} \left( Y^{p-1} \mu_{d} ; ({\bm a}, {\bm j})\right)$ 
in $\mathcal{U}_{r+1}$. It follows from the 
first claim in the step that each of the elements 
$$\mu_{\pm \sum_{i=0}^{r-1} p^i \left( \iota(-a_i,j_i)+2n^{(0)}(-a_i,j_i) \right)
-p^r d}^{(r+1)}X^{(p^r)p-1}
=\mu_{\pm n'-p^r d}^{(r+1)}X^{(p^r)p-1}$$ 
is a factor in a factorization  of 
$Z^{({\bm 0})} \left( \mu_{-d} X^{p-1} ; (-{\bm a}, {\bm j})\right)$ in $\mathcal{U}_{r+1}$ 
(recall from the proof of Proposition \ref{PropertiesOfB2} that $\iota(-a_i,j_i)=-b_i$). Since 
\begin{align*}
\mathcal{T}_2 \left( Z^{({\bm 0})} \left( \mu_{-d} X^{p-1} ; (-{\bm a}, {\bm j})\right) \right) 
&= \mathcal{T}_2 \left( 
\mu_{-\sum_{i=0}^{r-1}p^i b_i -p^r d}^{(r+1)} X^{(p^r)p-1} 
B^{({\bm 0})}(-{\bm a}, {\bm j})\right) \\
&= B^{({\bm 0})}({\bm a}, {\bm j})Y^{(p^r)p-1} 
\mu_{\sum_{i=0}^{r-1}p^i b_i +p^r d}^{(r+1)} \\
&= Y^{(p^r)p-1} 
\mu_{\sum_{i=0}^{r-1}p^i b_i +p^r d}^{(r+1)} B^{({\bm 0})}({\bm a}, {\bm j}) \\
&= Z^{({\bm 0})} \left( Y^{p-1}\mu_{d}  ; ({\bm a}, {\bm j})\right)
\end{align*}
by Propositions  \ref{PropertiesOfBZ} (iv), (v), \ref{PropertiesOfB2} (i), (ii), 
and \ref{CommProp} and since 
$$\mathcal{T}_2 \left( \mu_{\pm n'-p^r d}^{(r+1)}X^{(p^r)p-1} \right)
= Y^{(p^r)p-1} \mu_{\mp n'+p^r d}^{(r+1)}, 
$$
the second claim follows. \\

\noindent {\bf Step 2.} For an integer $\lambda$ with $0 \leq \lambda \leq \nu-1$, set 
$$z_{+,\lambda}=\mu_{\sum_{i=0}^{r-1} p^i\left(b_i + 2n^{(0)}(a_i,j_i)\right)+p^r a' -2p^r\lambda }^{(r+1)},$$
$$z_{-,\lambda}=\mu_{-\sum_{i=0}^{r-1} p^i\left(b_i + 2n^{(0)}(a_i,j_i)\right)+p^r a' -2p^r\lambda }^{(r+1)}.$$ 
Then there is an integer $\lambda$  with $0 \leq \lambda \leq \nu-1$ 
such that one of the two elements 
$z_{+,\lambda} X^{(p^r)p-\nu}$ and $z_{-,\lambda} X^{(p^r)p-\nu}$
lies in the two-sided ideal 
$\mathcal{U}_{r+1} h(\nu,r) X^{(p^r)p-\nu}\mathcal{U}_{r+1}$
of $\mathcal{U}_{r+1}$ 
and  there is an integer $\lambda$  with $0 \leq \lambda \leq \nu-1$ 
such that one of the two elements 
$Y^{(p^r)p-\nu} z_{+,\lambda}$ and $Y^{(p^r)p-\nu} z_{-,\lambda}$ 
lies in the two-sided ideal 
$\mathcal{U}_{r+1} Y^{(p^r)p-\nu}h(\nu,r)\mathcal{U}_{r+1}$ 
of $\mathcal{U}_{r+1}$. \\

Set $n'=\sum_{i=0}^{r-1}p^i \left( b_i+2n^{(0)}(a_i,j_i)\right)$  as in the proof of Step 1. 
Since  $z_{+, \lambda}$ and $z_{-, \lambda }$ are $\mathcal{U}_{r+1}^0$-weight vectors, there  
are scalars 
$c_{+, \lambda}, c_{-,\lambda} \in \mathbb{F}_p$ such that 
$h(\nu,r) z_{+,\lambda} =z_{+,\lambda} h(\nu,r) = 
c_{+,\lambda} z_{+,\lambda}$
and 
$h(\nu,r) z_{-,\lambda} =z_{-,\lambda} h(\nu,r) =
c_{-,\lambda} z_{-, \lambda}$.  
Note that 
$$c_{+,\lambda}={n' +p^r(a'+2\nu-2\lambda) \choose 2p^r \nu} +
{n' +p^r(a'+2\nu-2\lambda)-1 \choose 2p^r \nu},$$
$$c_{-,\lambda}={-n' +p^r(a'+2\nu-2\lambda) \choose 2p^r \nu} +
{-n' +p^r(a'+2\nu-2\lambda)-1 \choose 2p^r \nu}$$
and that $0 \leq n' \leq p^r-1$ by the remark of Definition \ref{Def2}. 
Therefore, we see that 
$$c_{+,\lambda}=2 {a'+2\nu-2\lambda \choose 2\nu},\ \ \ 
c_{-,\lambda}=2 {a'+2\nu-2\lambda-1 \choose 2\nu}$$
if $n'>0$ and 
$$c_{+,\lambda}=
c_{-,\lambda}= {a'+2\nu-2\lambda \choose 2\nu}+{a'+2\nu-2\lambda-1 \choose 2\nu}$$
if $n'=0$. 
To prove the claim in the step, we only have to show that there exists 
an integer $\lambda$ with $0 \leq \lambda \leq \nu-1$ such that 
$c_{+,\lambda} \neq 0$ or $c_{-,\lambda} \neq 0$ (in $\mathbb{F}_p$). Since 
$a' \not\equiv -1\ ({\rm mod}\ p)$, the integer $(a'+2\nu)\ {\bf mod}\ p $ lies in 
$\{ 0, \dots, p-1\} \backslash \{ 2 \nu -1\}$. 

Suppose that $2\nu \leq (a'+2\nu)\ {\bf mod}\ p \leq p-1$. Then if we take $\lambda=0$, 
we have $c_{+,\lambda}=2 {a'+2\nu \choose 2\nu} \neq 0$ in $\mathbb{F}_p$ 
when $n'>0$ and $c_{+,\lambda}= {a'+2\nu \choose 2\nu} + {a'+2\nu-1 \choose 2\nu}\neq 0$ in $\mathbb{F}_p$ 
when $n'=0$ by Lemma \ref{Lemma2}. 

Suppose that $0 \leq (a'+2\nu)\ {\bf mod}\ p \leq 2\nu-2$ and $n'>0$. Then if 
$(a'+2\nu)\ {\bf mod}\ p$ is even, taking 
$$\lambda=\dfrac{(a'+2\nu)\ {\bf mod}\ p}{2}\ (\leq \nu-1)$$
we have $c_{-,\lambda}=2  {a'+2\nu-2\lambda-1 \choose 2\nu}
=2 {p-1 \choose 2\nu} \neq 0$ in $\mathbb{F}_p$. If $(a'+2\nu)\ {\bf mod}\ p$ is odd, taking 
$$\lambda=\dfrac{\left( (a'+2\nu)\ {\bf mod}\ p\right)+1}{2}\ (\leq \nu-1)$$
we have $c_{+,\lambda}=2  {a'+2\nu-2\lambda \choose 2\nu}
=2 {p-1 \choose 2\nu} \neq 0$ in $\mathbb{F}_p$. 

Suppose that $0 \leq (a'+2\nu)\ {\bf mod}\ p \leq 2\nu-2$ and $n'=0$. Then if 
$(a'+2\nu)\ {\bf mod}\ p$ is even, taking 
$$\lambda=\dfrac{(a'+2\nu)\ {\bf mod}\ p}{2}\ (\leq \nu-1)$$
we have $c_{+,\lambda}=  {a'+2\nu-2\lambda \choose 2\nu} +
{a'+2\nu-2\lambda-1 \choose 2\nu}
= 0+{p-1 \choose 2\nu} \neq 0$ in $\mathbb{F}_p$. If $(a'+2\nu)\ {\bf mod}\ p$ is odd, taking 
$$\lambda=\dfrac{\left( (a'+2\nu)\ {\bf mod}\ p\right)+1}{2}\ (\leq \nu-1)$$
we have $c_{+,\lambda}=  {a'+2\nu-2\lambda \choose 2\nu} +
{a'+2\nu-2\lambda-1 \choose 2\nu} =
{p-1 \choose 2\nu} + {p-2 \choose 2\nu} \neq 0$ in $\mathbb{F}_p$ 
by Lemma \ref{Lemma2}. Therefore, Step 2 is proved. 

Now Steps 1 and 2 prove (i) and the proof is complete. $\square$ \\

\begin{Prop}\label{MainProp1}
Let $({\bm a}, {\bm j})= \left( (a_i,j_i)\right)_{i=0}^{r-1} \in \mathcal{P}_{\mathbb{Z}}^r$ and 
${\bm \varepsilon} =(\varepsilon_0, \dots, \varepsilon_{r-1}) \in 
\mathcal{X}_r({\bm a}, {\bm j})$. Suppose that there exists an integer 
$s \in \{ 0, \dots, r-1\}$ such that $\varepsilon_s=1$. Then the following hold. \\ \\
{\rm (i)} If $p$ is odd, the element $B^{({\bm \varepsilon})}({\bm a}, {\bm j})$ lies in the 
two-sided ideal  
$\mathcal{U}_r h(\nu,s) 
X^{(p^s)p-\nu} \mathcal{U}_r$ of 
$\mathcal{U}_r$  
for any $\nu \in \{ 1,2,\dots, (p-1)/2\}$. \\ \\
{\rm (ii)} If $p=2$,  the element $B^{({\bm \varepsilon})}({\bm a}, {\bm j})$ lies in the 
two-sided ideal 
$\mathcal{U}_r \mu_{m}^{(s+1)} X^{(m)} X^{(2^s)}\mathcal{U}_r$ of 
$\mathcal{U}_r$, where 
$m= \sum_{l=0}^{s-1}2^l(a_l\ {\bf mod}\ 2) $. 
\end{Prop}

\noindent {\itshape Proof.} 
The assumptions ${\bm \varepsilon} \in \mathcal{X}_r({\bm a}, {\bm j})$ and  
$\varepsilon_s=1$ imply that 
$(a_s,j_s)$ does not satisfy {\rm (E)}. 
If $r=1$ (and $s=0$), the proposition follows from Lemma 
\ref{Lemma4} (note that if $p=2$, then 
$(a_0\ {\bf mod}\ 2, j_0)=(0,1/2)$, $n^{(0)}(a_0,j_0)=0$, and $m=0$). Assume  
that $r \geq 2$. Since 
$$B^{({\bm \varepsilon})}({\bm a}, {\bm j})=B^{(\varepsilon_0)}(a_0,j_0) z$$
for some nonzero element $z \in \mathcal{U}_r$ by Proposition \ref{PropertiesOfZ} (ii), 
the result for $s=0$ follows from that for $r=1$. So we may 
assume that $s \geq 1$. 

Suppose first that $p=2$. 
Set ${\bm \varepsilon}'=(\varepsilon_1, \dots, \varepsilon_{r-1})$ and  
$({\bm a}', {\bm j}')= \left( (a_i, j_i)\right)_{i=1}^{r-1} $. Then we have 
$$B^{({\bm \varepsilon})}({\bm a}, {\bm j})= B^{(\varepsilon_0)}(a_0,j_0) Z^{(0)}
\left( B^{({\bm \varepsilon}')}({\bm a}', {\bm j}'); (a_0,j_0)\right).$$
By induction on $r$, the element 
$B^{({\bm \varepsilon}')}({\bm a}', {\bm j}')$ lies in 
$$\mathcal{U}_{r-1} \mu_{m'}^{(s)} X^{(m')} X^{(2^{s-1})}\mathcal{U}_{r-1},$$ where 
$m'= \sum_{l=0}^{s-2}2^l(a_{l+1}\ {\bf mod}\ 2) $. On the other hand, 
since $(a_0\ {\bf mod}\ 2)+2m'=m$, 
we see that 
$Z^{(0)} \left( \mu_{m'}^{(s)} X^{(m')} X^{(2^{s-1})} ; (a_0,j_0)\right)$ lies in 
$\mathcal{U}_{s+1} \mu_m^{(s+1)} X^{(m)} X^{(2^s)} \mathcal{U}_{s+1}$ by 
Lemma \ref{Lemma7} (ii). 
Therefore, by Proposition \ref{PropertiesOfZ} (iv) the element 
$B^{({\bm \varepsilon})}({\bm a}, {\bm j})$ 
must lie in 
$\mathcal{U}_{r} \mu_m^{(s+1)} X^{(m)} X^{(2^s)} \mathcal{U}_{r}$ and 
(ii) is proved.

From now on, assume that $p$ is odd. Fix $\nu \in \{ 1,2, \dots, (p-1)/2\}$. 
Set ${\bm \varepsilon}''=(\varepsilon_0, \dots, \varepsilon_{s-1})$, 
${\bm \varepsilon}'''=(\varepsilon_s, \dots, \varepsilon_{r-1})$, 
$({\bm a}'', {\bm j}'')= \left( (a_i, j_i)\right)_{i=0}^{s-1} $, and 
$({\bm a}''', {\bm j}''')= \left( (a_i, j_i)\right)_{i=s}^{r-1} $. By Proposition \ref{PropertiesOfZ} (ii) 
there is a nonzero element $z' \in \mathcal{U}_{r-s}$ such that 
$$B^{({\bm \varepsilon}''')}({\bm a}''', {\bm j}''')=B^{(1)}(a_s,j_s) z'.$$
Note from Lemma \ref{Lemma4} that $B^{(1)}(a_s,j_s)$ lies in 
$\mathcal{U}_1\mu_{a_s+2n^{(0)}(a_s,j_s)} X^{p-1} \mathcal{U}_1$. Since 
$(a_s,j_s)$ does not satisfy (E) (i.e. $j_s \neq 0$), we have  
$a_s+2n^{(0)}(a_s,j_s) \not\equiv -1\ ({\rm mod}\ p)$. Hence by 
Proposition \ref{PropertiesOfBZ} (vii) and Lemma 
\ref{Lemma7} (i)  the element 
$Z^{({\bm 0}_s)} \left( B^{(1)}(a_s,j_s) ; ({\bm a}'', {\bm j}'')\right)$ must lie in 
$\mathcal{U}_{s+1} h(\nu,s)X^{(p^s)p-\nu} \mathcal{U}_{s+1}$. 
Therefore, the element 
\begin{align*}
B^{({\bm \varepsilon})}({\bm a}, {\bm j})
&= B^{({\bm \varepsilon}'')}({\bm a}'', {\bm j}'')
Z^{({\bm 0}_s)} \left( B^{({\bm \varepsilon}''')}({\bm a}''', {\bm j}''') ; 
({\bm a}'', {\bm j}'') \right) \\
&= B^{({\bm \varepsilon}'')}({\bm a}'', {\bm j}'')
Z^{({\bm 0}_s)} \left( B^{(1)}(a_s,j_s) ; ({\bm a}'', {\bm j}'') \right)
Z^{({\bm 0}_s)} \left( z' ; ({\bm a}'', {\bm j}'') \right)
\end{align*} 
lies in $\mathcal{U}_{r} h(\nu,s)X^{(p^s)p-\nu} \mathcal{U}_{r}$ and (i) 
follows.  $\square$ \\

\begin{Prop}\label{MainProp2}
Let $({\bm a}, {\bm j})= \left( (a_i,j_i)\right)_{i=0}^{r-1} \in \mathcal{P}_{\mathbb{Z}}^r$ and 
${\bm \varepsilon} =(\varepsilon_0, \dots, \varepsilon_{r-1}) \in 
\mathcal{X}_r({\bm a}, {\bm j})$. Let  
${\bm t}=(t_0,\dots, t_{r-1})$ be an element of $\mathbb{Z}^r$ satisfying 
$-\widetilde{n}^{(\varepsilon_i+1)}(a_i,j_i) \leq t_i \leq n^{(\varepsilon_i+1)}(a_i,j_i)$ 
for any $i \in \{ 0, \dots, r-1\}$ (then 
$B^{({\bm \varepsilon})}\left( ({\bm a}, {\bm j}); {\bm t}\right) \neq 0$). 
Suppose that there is an integer $s \in \{0,\dots, r-1 \}$ 
such that $\varepsilon_s=0$ and 
$$n^{(0)}(a_s,j_s) < t_s \leq n^{(1)}(a_s,j_s)\ \  \mbox{or}\ \ 
-\widetilde{n}^{(1)}(a_s,j_s) \leq t_s < -\widetilde{n}^{(0)}(a_s,j_s)$$
(these occur only if $(a_s,j_s)$ does not satisfy {\rm (E)}). Then the following hold. \\ 

\noindent{\rm (i)} Suppose that $p$ is odd. For any $\nu \in \{ 1,2,\dots, (p-1)/2\}$, the element 
$B^{({\bm \varepsilon})}\left( ({\bm a}, {\bm j}); {\bm t}\right)$ lies in the two-sided ideal   
$$
\left\{ \begin{array}{ll}
{\mathcal{U}_r h(\nu,s) X^{(p^s)p-\nu} \mathcal{U}_r} & 
\mbox{if $n^{(0)}(a_s,j_s) < t_s \leq n^{(1)}(a_s,j_s)$,} \\
{\mathcal{U}_r Y^{(p^s)p-\nu} h(\nu,s) \mathcal{U}_r} &
\mbox{if $-\widetilde{n}^{(1)}(a_s,j_s) \leq t_s < -\widetilde{n}^{(0)}(a_s,j_s)$}
\end{array} \right.
$$
of $\mathcal{U}_r$.\\ 

\noindent {\rm (ii)} Suppose that  $p=2$. Then the element 
$B^{({\bm \varepsilon})}\left( ({\bm a}, {\bm j}); {\bm t}\right) $ 
lies in the two-sided ideal 
$$
\left\{ \begin{array}{ll}
{\mathcal{U}_r \mu_{m}^{(s+1)} X^{(m)} X^{(2^s)}\mathcal{U}_r} & 
\mbox{if $n^{(0)}(a_s,j_s) < t_s \leq n^{(1)}(a_s,j_s)$} \\
   &  \mbox{(i.e. if $(a_s\ {\bf mod}\ 2, j_s)=(0,1/2)$ and $t_s=1$),} \\
{\mathcal{U}_r Y^{(2^s)}Y^{(m)}\mu_{m}^{(s+1)}  \mathcal{U}_r} &
\mbox{if $-\widetilde{n}^{(1)}(a_s,j_s) \leq t_s < -\widetilde{n}^{(0)}(a_s,j_s)$} \\
   &  \mbox{(i.e. if $(a_s\ {\bf mod}\ 2, j_s)=(0,1/2)$ and $t_s=-1$)}
\end{array} \right.
$$
of $\mathcal{U}_r$, where $m= \sum_{l=0}^{s-1}2^l(a_l\ {\bf mod}\ 2) $. 
\end{Prop}

\noindent {\itshape Proof.} For $r=1$ (and $s=0$), (i) and (ii) have been proved 
in Lemma \ref{Lemma5}. 
Assume that $r \geq 2$. Note that 
$$B^{({\bm \varepsilon})}\left( ({\bm a}, {\bm j}); {\bm t}\right)
= B^{(\varepsilon_0)}\left( (a_0,j_0); t_0\right) Z^{(0)} 
\left(B^{({\bm \varepsilon}')}\left( ({\bm a}', {\bm j}'); {\bm t}' \right);(a_0,j_0) \right),$$
where $ ({\bm a}', {\bm j}')= \left( (a_i,j_i)\right)_{i=1}^{r-1}$, 
${\bm \varepsilon}'=(\varepsilon_1, \dots, \varepsilon_{r-1})$, and 
${\bm t}'=(t_1, \dots, t_{r-1})$. Then Lemma \ref{Lemma5} shows the claims for 
$s=0$ in (i) and (ii). 
 Therefore, we may assume that $s \geq 1$. 

Suppose that $p=2$.  Since $B^{({\bm \varepsilon}')}\left( ({\bm a}', {\bm j}'); {\bm t}' \right)$ 
lies in $\mathcal{U}_{r-1} \mu_{m'}^{(s)}X^{(m')} X^{(2^{s-1})} \mathcal{U}_{r-1}$ if  
$n^{(0)}(a_s,j_s) < t_s \leq n^{(1)}(a_s,j_s)$ and in  $\mathcal{U}_{r-1} Y^{(2^{s-1})}Y^{(m')}\mu_{m'}^{(s)}  \mathcal{U}_{r-1}$ if 
$-\widetilde{n}^{(1)}(a_s,j_s) \leq t_s < -\widetilde{n}^{(0)}(a_s,j_s)$ 
by induction on $r$, where $m'= \sum_{l=0}^{s-2} 2^l(a_{l+1}\ {\bf mod}\ 2)$, 
Proposition \ref{PropertiesOfZ} (iv) and Lemma \ref{Lemma7} (ii) imply that the element 
$B^{({\bm \varepsilon})}\left( ({\bm a}, {\bm j}); {\bm t}\right)$ 
must lie in $\mathcal{U}_{r} \mu_{m}^{(s+1)}X^{(m)} X^{(2^{s})} \mathcal{U}_{r}$ 
if $n^{(0)}(a_s,j_s) < t_s \leq n^{(1)}(a_s,j_s)$ and in 
$\mathcal{U}_{r} Y^{(2^{s})}Y^{(m)}\mu_{m}^{(s+1)}  \mathcal{U}_{r}$ 
if $-\widetilde{n}^{(1)}(a_s,j_s) \leq t_s < -\widetilde{n}^{(0)}(a_s,j_s)$, and (ii) is proved. 

From now on, assume that $p$ is odd. Fix 
$\nu \in \{ 1,2,\dots, (p-1)/2\}$. Set 
${\bm \varepsilon}''=(\varepsilon_0, \dots, \varepsilon_{s-1})$, 
${\bm \varepsilon}'''=(\varepsilon_s, \dots, \varepsilon_{r-1})$, 
$({\bm a}'', {\bm j}'')= \left( (a_i, j_i)\right)_{i=0}^{s-1} $, 
$({\bm a}''', {\bm j}''')= \left( (a_i, j_i)\right)_{i=s}^{r-1} $, 
${\bm t}''=(t_0, \dots, t_{s-1})$, and 
${\bm t}'''=(t_s, \dots, t_{r-1})$.  Then 
there is a nonzero element $z' \in \mathcal{U}_{r-s}$ such that 
$$B^{({\bm \varepsilon}''')}\left( ({\bm a}''', {\bm j}'''); {\bm t}''' \right)
=B^{(0)}\left( (a_s,j_s); t_s \right) z'.$$
Indeed, we can take $z'=1$ if $s=r-1$ and 
$$
z'=Z^{(0)} 
\left(B^{(\varepsilon_{s+1}, \dots, \varepsilon_{r-1})}\left( ((a_i,j_i))_{i=s+1}^{r-1}; 
(t_{s+1}, \dots, t_{r-1}) \right);(a_s,j_s) \right)
$$
if $s<r-1$. 
Note from Lemma \ref{Lemma5} that 
$B^{(0)}\left( (a_s,j_s);t_s\right)$ lies in 
$\mathcal{U}_1 \mu_{a_s+2n^{(1)}(a_s,j_s)} X^{p-1} \mathcal{U}_1$ 
if $n^{(0)}(a_s,j_s) < t_s \leq n^{(1)}(a_s,j_s)$ and in 
$\mathcal{U}_1  Y^{p-1} \mu_{a_s+2n^{(0)}(a_s,j_s)}\mathcal{U}_1$ 
if $-\widetilde{n}^{(1)}(a_s,j_s) \leq t_s < -\widetilde{n}^{(0)}(a_s,j_s)$.
Since $(a_s,j_s)$ does not satisfy (E) 
(i.e. $j_s \neq 0$),  neither $a_s+2n^{(1)}(a_s,j_s)$ nor $a_s+2n^{(0)}(a_s,j_s)$ 
is  congruent to $-1$ modulo 
$p$. Hence by Proposition \ref{PropertiesOfBZ} (vii) and Lemma \ref{Lemma7} (i) the element 
$Z^{({\bm 0}_s)} \left( B^{(0)}\left( (a_s,j_s);t_s \right) ; ({\bm a}'', {\bm j}'')\right)$ must lie in 
$\mathcal{U}_{s+1} h(\nu,s)X^{(p^s)p-\nu} \mathcal{U}_{s+1}$ if 
$n^{(0)}(a_s,j_s) < t_s \leq n^{(1)}(a_s,j_s)$ and in 
$\mathcal{U}_{s+1} Y^{(p^s)p-\nu}h(\nu,s) \mathcal{U}_{s+1}$ 
if  $-\widetilde{n}^{(1)}(a_s,j_s) \leq t_s < -\widetilde{n}^{(0)}(a_s,j_s)$. 
Therefore, the element 
$$B^{({\bm \varepsilon})}(({\bm a}, {\bm j});{\bm t})
= B^{({\bm \varepsilon}'')}\left( ({\bm a}'', {\bm j}'') ; {\bm t}'' \right)
Z^{({\bm 0}_s)} \left( B^{(0)}\left( (a_s,j_s);t_s \right) ; ({\bm a}'', {\bm j}'') \right)
Z^{({\bm 0}_s)} \left( z' ; ({\bm a}'', {\bm j}'') \right)$$
lies in $\mathcal{U}_{r} h(\nu,s)X^{(p^s)p-\nu} \mathcal{U}_{r}$ 
if $n^{(0)}(a_s,j_s) < t_s \leq n^{(1)}(a_s,j_s)$ and in 
$\mathcal{U}_{r} Y^{(p^s)p-\nu}h(\nu,s) \mathcal{U}_{r}$ 
if $-\widetilde{n}^{(1)}(a_s,j_s) \leq t_s < -\widetilde{n}^{(0)}(a_s,j_s)$, and (i) 
follows.  $\square$  \\ 

Let us now turn to the proof of the main theorem. \\

\noindent {\itshape Proof of Theorem \ref{MainThm}.}  
Let $\mathcal{I}$ be the two-sided ideal of 
$\mathcal{U}_r$ generated by the set in (i) or (ii). 
Fix the elements $\nu_l$ with $l \in \{1,2,\dots, (p-1)/2\}$ in 
(i) if  $p$ is odd. We proceed in steps. \\ 

\noindent {\bf Step 1.} $\mathcal{I} \subseteq {\rm rad} \mathcal{U}_r$. \\ 

Since it has been proved for $r=1$ in Lemma \ref{Lemma6}, we may assume that 
$r \geq 2$. 

We know  that  the set 
$$\left\{ B^{({\bm 1})}\left( ({\bm a}, {\bm j}); {\bm t} \right) \left| 
\begin{array}{l}
\mbox{$({\bm a}, {\bm j}) =\left( (a_i, j_i)\right)_{i=0}^{r-1}\in \mathcal{P}^r$,\ 
${\bm t}=(t_0, \dots, t_{r-1})$, } \\
\mbox{$-\widetilde{n}^{(0)}(a_i,j_i) \leq t_i \leq n^{(0)}(a_i,j_i)$, $\forall i$}
\end{array}
\right.\right\}$$
forms a $k$-basis of the $\mathcal{U}_r$-module ${\rm soc}_{\mathcal{U}_r} \mathcal{U}_r$
(see Theorem 5.1 and the remark of Theorem 5.3 in \cite{yoshii22}). 
Choose a basis element $B^{({\bm 1})}\left( ({\bm a}, {\bm j}); {\bm t} \right)$ 
arbitrarily, where ${\bm t}=(t_0, \dots, t_{r-1})$ and 
$-\widetilde{n}^{(0)}(a_i,j_i) \leq t_i \leq n^{(0)}(a_i,j_i)$ for any $i$. 
It is enough to check that 
the element $B^{({\bm 1})}\left( ({\bm a}, {\bm j}); {\bm t} \right)$  
is annihilated by 
any element in the sets generating $\mathcal{I}$. Note that 
there exists a nonzero 
element $z \in \mathcal{U}_r$ such that 
$$B^{({\bm 1})}\left( ({\bm a}, {\bm j}); {\bm t} \right)=
B^{(1)}\left( (a_0,j_0); t_0\right) z.$$
Indeed, we can take 
$z=Z^{(0)} 
\left(B^{({\bm 1}_{r-1})}\left( ({\bm a}', {\bm j}'); {\bm t}' \right);(a_0,j_0) \right)$, 
where $ ({\bm a}', {\bm j}')= \left( (a_i,j_i)\right)_{i=1}^{r-1}$ and 
${\bm t}'=(t_1, \dots, t_{r-1})$. 
Since $B^{(1)}\left( (a_0,j_0); t_0\right)$ 
lies in ${\rm soc}_{\mathcal{U}_1} \mathcal{U}_1$, Lemma \ref{Lemma6} shows that 
$$h(\nu_0,0) X^{p-\nu_0} 
B^{({\bm 1})}\left( ({\bm a}, {\bm j}); {\bm t} \right) = 
Y^{p-\nu_0} h(\nu_0,0) 
B^{({\bm 1})}\left( ({\bm a}, {\bm j}); {\bm t} \right)=0$$
if $p$ is odd and 
$$\mu_0 X B^{({\bm 1})}\left( ({\bm a}, {\bm j}); {\bm t} \right)= 
Y \mu_0  B^{({\bm 1})}\left( ({\bm a}, {\bm j}); {\bm t} \right)=0$$
if $p=2$. Assume that $1 \leq i \leq r-1$. 

Suppose that $p=2$ and consider the element 
$\mu_m^{(i+1)} X^{(m)} X^{(2^i)}B^{({\bm 1})}\left( ({\bm a}, {\bm j}); {\bm t} \right)$ 
for an integer $m$ with $0 \leq m \leq 2^i-1$.  
Let $m= \sum_{l=0}^{i-1} 2^l m_l $ with $m_l \in \{0,1 \}$ be the $2$-adic expansion of $m$. 
By \cite[Lemma 4.7]{yoshii22}, 
we know that $X^{(m)} X^{(2^i)} B^{({\bm 1})}\left( ({\bm a}, {\bm j}); {\bm t} \right)$ is 
a nonzero scalar multiple of 
$$B^{({\bm 1})}\left( ({\bm a}, {\bm j}); {\bm t}+{\bm e}_{i+1} + 
\sum_{l=0}^{i-1} m_l {\bm e}_{l+1} \right).$$ 
If the element $B^{({\bm 1})}\left( ({\bm a}, {\bm j}); {\bm t}+{\bm e}_{i+1} + 
\sum_{l=0}^{i-1} m_l {\bm e}_{l+1} \right)$ is not zero, not only 
it has $\mathcal{U}_{i+1}^{0}$-weight 
$$\sum_{l=0}^i 2^l b_l +2 \sum_{l=0}^i 2^l t_l + 2 \cdot 2^i+2m$$
but also we must have   
$t_i +1 \leq n^{(0)}(a_i,j_i)$ and $t_l +m_l \leq n^{(0)}(a_l, j_l)$ for $l \in \{0, \dots, i-1\}$  
(see \cite[Proposition 4.5]{yoshii22}). 
Now we shall show that 
$$\mbox{$b_i= \pm 1$ and $-1 \leq b_l+2t_l+m_l \leq 1$ for any 
$l \in \{ 0, \dots, i-1\}$.} $$
Suppose that $(a_i\ {\bf mod}\ 2, j_i)=(0,1/2)$. Then we have 
$n^{(0)}(a_i,j_i)=\widetilde{n}^{(0)}(a_i,j_i)=0$ and $t_i=0$, which does not satisfy 
$t_i+1 \leq n^{(0)}(a_i,j_i)$. So $(a_i\ {\bf mod}\ 2, j_i)$ must be $(1,0)$ or $(1,1)$ and hence 
we have $b_i = \pm 1$. In turn, let $l \in \{0, \dots, i-1\}$. If $m_l=0$, we have 
$$-1 \leq b_l -2 \widetilde{n}^{(0)}(a_l,j_l) \leq b_l + 2t_l +m_l \leq 
b_l +2 n^{(0)}(a_l,j_l) \leq 1.$$
On the other hand, if $m_l=1$, 
$(a_l\ {\bf mod}\ 2, j_l)$ must be $(1,0)$ or $(1,1)$ by the same argument as above.  
Then  $(b_l, t_l)$ must be $(-1,0)$ or $(1,-1)$ and hence we have 
$ b_l + 2t_l +m_l =0$. Therefore, the above claim follows. 
 
Then  we have 
$$-2^i+1 \leq \sum_{l=0}^{i-1}2^l(b_l+2t_l+m_l) \leq 2^i-1$$
and hence
$$
\sum_{l=0}^i 2^l b_l +2 \sum_{l=0}^i 2^l t_l + 2 \cdot 2^i+2m
\equiv \sum_{l=0}^{i-1}2^l(b_l+2t_l+m_l)+2^i +m 
\not\equiv m\ ({\rm mod}\ 2^{i+1}).
$$
Therefore, the $\mathcal{U}_{i+1}^{0}$-weight of nonzero 
$B^{({\bm 1})}\left( ({\bm a}, {\bm j}); {\bm t}+{\bm e}_{i+1} + 
\sum_{l=0}^{i-1} m_l {\bm e}_{l+1} \right)$ is not congruent to $m$ modulo $2^{i+1}$ and 
hence we obtain 
$$\mu_m^{(i+1)}X^{(m)}X^{(2^i)}B^{({\bm 1})}\left( ({\bm a}, {\bm j}); {\bm t} \right)=0.$$

We also note from this fact that 
$$\mu_m^{(i+1)}X^{(m)}X^{(2^i)}B^{({\bm 1})}\left( (-{\bm a}, {\bm j}); -{\bm t} \right)=0$$
since $-\widetilde{n}^{(0)}(-a_l,j_l) \leq -t_l \leq n^{(0)}(-a_l,j_l)$ for all 
$l \in \{ 0,\dots, r-1\}$. Since 
\begin{align*}
\lefteqn{\mathcal{T}_1 \left( \mu_m^{(i+1)}X^{(m)}X^{(2^i)}B^{({\bm 1})}
\left( (-{\bm a}, {\bm j}); -{\bm t} \right) \right)} \\
&= (-1)^{m+2^i +\sum_{l=0}^{r-1} 2^l |t_l| } 
\mu_{-m}^{(i+1)}Y^{(m)}Y^{(2^i)}B^{({\bm 1})}
\left( ({\bm a}, {\bm j}); {\bm t} \right) \\
&= (-1)^{m+2^i +\sum_{l=0}^{r-1} 2^l |t_l| } 
Y^{(m)}Y^{(2^i)}\mu_{m}^{(i+1)}B^{({\bm 1})}
\left( ({\bm a}, {\bm j}); {\bm t} \right) 
\end{align*}
by Propositions \ref{PropertiesOfB2} (ii) and \ref{FormulasOfMu} (i), we obtain 
$$Y^{(m)}Y^{(2^i)}\mu_{m}^{(i+1)}B^{({\bm 1})}
\left( ({\bm a}, {\bm j}); {\bm t} \right)=0.$$
Therefore, 
we have proved that $\mathcal{I} \subseteq {\rm rad} \mathcal{U}_r$ for $p=2$. 

Suppose that $p$ is odd. To prove 
$\mathcal{I} \subseteq {\rm rad} \mathcal{U}_r$, it remains  to show that 
$$h(\nu_i,i)X^{(p^i)p-\nu_i} B^{({\bm 1})}\left( ({\bm a}, {\bm j}); {\bm t} \right)=
Y^{(p^i)p-\nu_i} h(\nu_i,i) B^{({\bm 1})}\left( ({\bm a}, {\bm j}); {\bm t} \right)=0$$
for $1 \leq i \leq r-1$. It is enough to show that  
$${H+2p^i\nu_i -\delta \choose 2p^i\nu_i }X^{(p^i)p-\nu_i} 
B^{({\bm 1})}\left( ({\bm a}, {\bm j}); {\bm t} \right)=
Y^{(p^i)p-\nu_i} {H+2p^i\nu_i -\delta \choose 2p^i\nu_i }
B^{({\bm 1})}\left( ({\bm a}, {\bm j}); {\bm t} \right)=0
$$ 
for $\delta \in \{0,1 \}$ and $1 \leq i \leq r-1$. 
We know that $X^{(p^i)p-\nu_i} B^{({\bm 1})}\left( ({\bm a}, {\bm j}); {\bm t} \right)$ is 
a nonzero scalar multiple of 
$$B^{({\bm 1})}\left( ({\bm a}, {\bm j}); {\bm t} +(p-\nu_i){\bm e}_{i+1}\right)$$
by \cite[Lemma 4.7]{yoshii22}.  
If the element  $B^{({\bm 1})}\left( ({\bm a}, {\bm j}); {\bm t} +(p-\nu_i){\bm e}_{i+1}\right)$ 
is not zero, we have not only 
\begin{align*}
\lefteqn{{H+2p^i\nu_i -\delta \choose 2p^i\nu_i }
B^{({\bm 1})}\left( ({\bm a}, {\bm j}); {\bm t} +(p-\nu_i){\bm e}_{i+1}\right)} \\
&= {\sum_{l=0}^{i}p^l(b_l+2t_l) +p^{i+1}-\delta \choose 2p^i\nu_i }
B^{({\bm 1})}\left( ({\bm a}, {\bm j}); {\bm t} +(p-\nu_i){\bm e}_{i+1}\right) 
\end{align*}
but also $t_i+p -\nu_i \leq n^{(0)}(a_i,j_i)$ by \cite[Proposition 4.5]{yoshii22}. Since  
$$-p+1 \leq b_l-2\widetilde{n}^{(0)}(a_l,j_l) \leq b_l+2t_l \leq b_l+2n^{(0)}(a_l,j_l) \leq p-1$$
for  $0 \leq l \leq i-1$ and 
$$-p+1 \leq b_i-2\widetilde{n}^{(0)}(a_i,j_i) \leq b_i+2t_i \leq b_i+2n^{(0)}(a_i,j_i) 
+2(\nu_i-p)\leq -p+2\nu_i-1$$
by the remark of Definition \ref{Def2}, we see that 
$$1 \leq \sum_{l=0}^{i}p^l(b_l+2t_l) +p^{i+1} \leq 2p^i\nu_i  -1.$$
Thus we have ${\sum_{l=0}^{i}p^l(b_l+2t_l) +p^{i+1}-\delta \choose 2p^i\nu_i } = 0$ in 
$\mathbb{F}_p$. Therefore, we obtain 
$${H+2p^i\nu_i -\delta \choose 2p^i\nu_i }X^{(p^i)p-\nu_i} B^{({\bm 1})}
\left( ({\bm a}, {\bm j}); {\bm t} \right) =0.$$
We also note from this result and 
$-\widetilde{n}^{(0)}(-a_l,j_l) \leq -t_l \leq n^{(0)}(-a_l, j_l)$ for all 
$l \in \{0,\dots, r-1 \}$ that 
$${H+2p^i\nu_i -\delta \choose 2p^i\nu_i }X^{(p^i)p-\nu_i} 
B^{({\bm 1})}\left( (-{\bm a}, {\bm j}); -{\bm t} \right) =0.$$
Since it is easy to check that 
$\mathcal{T}_1 \left( {H+2p^i\nu_i -\delta \choose 2p^i\nu_i } \right)
= {H-(1-\delta) \choose 2p^i\nu_i } $, we obtain 
\begin{align*}
\lefteqn{\mathcal{T}_1 \left( {H+2p^i\nu_i -\delta \choose 2p^i\nu_i }
X^{(p^i)p-\nu_i} B^{({\bm 1})}\left( (-{\bm a}, {\bm j}); -{\bm t} \right) \right)} \\
&= (-1)^{p^i(p-\nu_i)+\sum_{l=0}^{r-1}p^l |t_l| } 
{H-(1-\delta) \choose 2p^i\nu_i } Y^{(p^i)p-\nu_i} 
B^{({\bm 1})}\left( ({\bm a}, {\bm j}); {\bm t} \right) \\
&= (-1)^{p^i(p-\nu_i)+\sum_{l=0}^{r-1}p^l |t_l| } 
Y^{(p^i)p-\nu_i} {H+2p^i\nu_i -(1-\delta) \choose 2p^i\nu_i }
B^{({\bm 1})}\left( ({\bm a}, {\bm j}); {\bm t} \right) 
\end{align*}
by Proposition \ref{PropertiesOfB2} (ii) and hence 
$$Y^{(p^i)p-\nu_i} {H+2p^i\nu_i -\delta \choose 2p^i\nu_i }
B^{({\bm 1})}\left( ({\bm a}, {\bm j}); {\bm t} \right)=0$$
for $\delta \in \{0,1\}$. 
Therefore, 
we have proved that $\mathcal{I} \subseteq {\rm rad} \mathcal{U}_r$ for odd $p$ 
and Step 1 is proved. \\ 

\noindent {\bf Step 2.} $\mathcal{I} \supseteq {\rm rad} \mathcal{U}_r$. \\

Set 
$$\Xi_r({\bm a},{\bm j})= \Theta_r \left( ({\bm a},{\bm j}), {\bm 0}\right) 
\backslash \left\{ ({\bm 0},{\bm t}({\bm 0})) \left| 
-\widetilde{n}^{(0)}(a_i,j_i) \leq t_i(0) \leq n^{(0)}(a_i,j_i),\ \forall i \right. \right\},$$
where ${\bm t}({\bm 0})=(t_0(0), \dots, t_{r-1}(0))$. 
Recall from Theorem \ref{BasisThm} that the elements 
$B^{({\bm \theta})}\left( ({\bm a}, {\bm j}); {\bm t}({\bm \theta})\right)$ 
with $({\bm a},{\bm j}) \in \mathcal{P}^r$ and 
$({\bm \theta}, {\bm t}({\bm \theta})) \in \Xi_r({\bm a},{\bm j})$ form a $k$-basis of 
${\rm rad} \mathcal{U}_r$. 
Then every such basis element 
$B^{({\bm \theta})}\left( ({\bm a}, {\bm j}); {\bm t}({\bm \theta})\right)$  
satisfies one of the following: 
\begin{itemize}
\item There is an integer $s \in \{0, \dots, r-1\}$ such that $\theta_s=0$ and 
$${n}^{(0)}(a_s,j_s) < t_s(\theta_s) \leq {n}^{(1)}(a_s,j_s)\ \mbox{ or }\  
-\widetilde{n}^{(1)}(a_s,j_s) \leq t_s(\theta_s) < -\widetilde{n}^{(0)}(a_s,j_s).$$ 
\item There is an integer $s \in \{0, \dots, r-1\}$ such that $\theta_s=1$. 
\end{itemize}  
(For the notation, see Definition \ref{Def2} (4).) 
Now Propositions \ref{MainProp1} and \ref{MainProp2} imply that in each case the element 
$B^{({\bm \theta})}\left( ({\bm a}, {\bm j}); {\bm t}({\bm \theta})\right)$ 
must lie in $\mathcal{I}$, and hence the claim follows. $\square$

\ \\
{\large {\bf Acknowledgment}} \\

The author would like to thank the referee for carefully reading the first version of this paper and giving  helpful  comments.

This work was supported by JSPS KAKENHI Grant Number JP18K03203.

\end{document}